\documentclass[11pt]{article}
\usepackage[T1]{fontenc}
\usepackage{amssymb,amsmath,amsfonts,dsfont,amsthm,bm}
\usepackage{float,rotfloat} 
\usepackage{mathrsfs}
\usepackage{multirow}
\usepackage{enumerate}
\usepackage{enumitem}
\usepackage{mathtools}
\usepackage{multicol}
\usepackage{hyperref,url}
\usepackage{graphicx}
\usepackage[hang,small,bf]{caption}
\usepackage[top=1in, bottom=1in, left=1.0in, right=1.0 in]{geometry}
\usepackage{longtable}
\usepackage{color,xcolor}
\usepackage[round,sort,comma]{natbib}
\newcounter{bibcount}
\usepackage[title]{appendix}
\usepackage{etoolbox} 
\usepackage[english]{babel}
\usepackage[autostyle, english = american]{csquotes}
\MakeOuterQuote{"}

\usepackage{cancel}

\patchcmd{\appendices}{\quad}{: }{}{} 

\newtheorem{thm}{Theorem}

\newtheorem{cor}{Corollary}

\newtheorem{note}{Note}
\newcommand{\ID}{\mathds{1}}
\newcommand{\lf}{\left\lfloor}
\newcommand{\rf}{\right\rfloor}
\newcommand{\ol}{\overline}
\newcommand{\lb}{\left[}
\newcommand{\rb}{\right]}
\newcommand{\lp}{\left(}
\newcommand{\rp}{\right)}

\newcommand{\R}{\mathbb{R}}

\definecolor{darkblue}{rgb}{0,0,0.55}

\definecolor{dkgreen}{rgb}{0,0.6,0}

\allowdisplaybreaks

\hypersetup{,colorlinks=true,allcolors=darkblue}

\makeatletter
\patchcmd{\@lbibitem}%
{\item[}%
{\item[\hfill\stepcounter{bibcount}{[\thebibcount]}}%
{}%
{}
\patchcmd{\@lbibitem}%
{\hfil \NAT@anchor {#2}{\NAT@num }]}%
{\hyper@natanchorstart{#2\@extra@b@citeb}\hyper@natanchorend]}
{}%
{}
\renewcommand\NAT@bibsetup[1]{%
\setlength{\leftmargin}{0.40in}
\setlength{\itemindent}{-\parindent}%
\setlength{\itemsep}{\bibsep}%
\setlength{\parsep}{\z@}}
\makeatother

\begin{document}   

\baselineskip 5mm

\thispagestyle{empty}

\begin{center}
{\Large
On the Asymptotic Normality of 
Trimmed and Winsorized $L$-statistics} 

\vspace{15mm}

{\large\sc
Chudamani Poudyal\footnote[1]{ 
Chudamani Poudyal, PhD, ASA, 
is an Assistant Professor 
in the Department of Mathematical Sciences, 
University of Wisconsin-Milwaukee, 
Milwaukee, WI 53211, USA. 
~~ {\em E-mail\/}: ~{\tt cpoudyal@uwm.edu}}}

\vspace{1mm}

{\large\em Department of Mathematical Sciences} \\[5pt]
{\large\em University of Wisconsin-Milwaukee}

\vspace{15mm}

\vspace{15mm}
\end{center}

\vspace{5mm}

\begin{quote}
{\bf\em Abstract\/}.

There are several ways to establish the asymptotic
normality of $L$-statistics,
which depend on the choice of the 
weights-generating function and 
the cumulative distribution selection 
of the underlying model. 
In this study, we focus on establishing 
computational formulas 
for the asymptotic variance 
of two robust $L$-estimators: 
the method of trimmed moments (MTM) and 
the method of winsorized moments (MWM). 
We demonstrate that two asymptotic approaches
for MTM are equivalent for a specific choice 
of the weights-generating function. 
These findings enhance the applicability 
of these estimators across various underlying distributions, 
making them effective tools in diverse statistical scenarios. 
Such scenarios include actuarial contexts,
such as payment-per-payment and payment-per-loss data scenarios, 
as well as in evaluating the asymptotic 
distributional properties of distortion risk measures.
The effectiveness of our methodologies depends on the 
availability of the cumulative distribution function, 
ensuring broad usability in various statistical environments.

\vspace{4mm}


\vspace{4mm}

\end{quote}

\newpage

\baselineskip 7mm
\setcounter{page}{2}

\section{Introduction}
\label{sec:Intro}
 
For a positive integer $n$, let $X_{1}, \ldots, X_{n}$
be an iid sample from a distribution with 
an unknown true underlying cumulative distribution function $F$, 
characterized by the parameter vector
$\bm{\theta} = (\theta_{1}, \cdots, \theta_{k})$. 
The primary motivation of most statistical inference efforts
is to estimate the parameter vector $\bm{\theta}$
from the available sample dataset. 
The corresponding order statistic 
of the sample are denoted by 
$X_{1:n} \le \ldots \le X_{n:n}$.
In order to estimate the 
parameter vector $\bm{\theta}$,
the statistics we are interested 
here is a {\em linear combination}
of the order statistics,
so the name $L$-statistics, 
in the form 
\begin{align}
\label{eqn:S1}
\widehat{\mu}_{j} 
& :=
\dfrac{1}{n}
\sum_{i=1}^{n}
J_{j} \left( \dfrac{i}{n+1} \right) 
h_{j}(X_{i:n}),
\quad 
1 \le j \le k,
\end{align}
where 
$J_{j} : [0,1] \to [0, \infty)$
represents a {\em weights-generating} 
function.
Both $J_{j}$ and $h_{j}$
are specially chosen functions,
\citep[see, e.g.,][]{MR4263275, pb22}
and are known that are specified
by the statistician.

The corresponding population 
quantities
are then given by 
\begin{align}
\label{eqn:P1}
\mu_{j}
& \equiv 
\mu_{j} 
\left( 
\bm{\theta}
\right) 
\equiv 
\mu_{j} 
\left( 
\theta_1, \ldots, \theta_k
\right) 
= 
\int_{0}^{1}
{J_{j}(u)H_{j}(u) \, du,
\quad 
1\leq{j}\leq{k}},
\mbox{ where } 
H_{j} := h_{j} \circ F^{-1}.
\end{align}

From Eq. \eqref{eqn:P1},
it clearly follows that
\begin{align}
\label{eqn:cgj7}
H_{j}'(u) du
&= 
dh_{j} \lp F^{-1}(u) \rp.
\end{align}

Like classical {\em method of moments}, 
$L$-estimators
\citep[see, e.g.,][]{MR1049304}
are found by matching sample 
$L$-moments, Eq. \eqref{eqn:S1}, 
with population $L$-moments, 
Eq. \eqref{eqn:P1},
for $j = 1, \ldots, k$,
and then solving the system of equations
with respect to $\theta_1, \ldots, \theta_k$.
The obtained solutions, 
which we denote by 
$
\widehat{\theta}_j 
=
\tau_j(\widehat{\mu}_{1}, \ldots, \widehat{\mu}_{k})$,
$1 \leq j \leq k$, are, by definition, 
the $L$-estimators of 
$\theta_1, \ldots, \theta_k$. 
Note that the functions $\tau_j$ are such that
$
\theta_j 
=
\tau_j(\mu_1(\boldsymbol{\theta}),
\ldots, \mu_k(\boldsymbol{\theta})).
$

For \(0 \le u \le 1\), 
we define \(\overline{u} = 1-u\),
and consider
\begin{align}
\label{eqn:AlphaFun1}
\alpha_{j}(u)
& = 
\frac{1}{\ol{u}}
\int_{u}^{1}
\ol{v} \ g_{j}(v) \, dv, 
\quad 
\mbox{where}
\quad 
g_{j}(v) 
:=
J_{j}(v)H_{j}'(v),
\quad 
1\leq{j}\leq{k}.
\end{align}

Furthermore, consider the estimators 
and true parameter vectors:
\[
\widehat{\bm\mu}
:=
\left( 
\widehat{\mu}_{1}, 
\widehat{\mu}_{1}, 
\ldots, 
\widehat{\mu}_{k}
\right)
\quad 
\mbox{and}
\quad 
\bm{\mu}
:=
\left(
\mu_{1}, 
\mu_{2}, 
\ldots, 
\mu_{k}
\right).
\]

Ideally, 
we expect that the sample vector 
$\widehat{\bm\mu}$
converges in distribution 
to the population vector $\bm{\mu}$.
As mentioned by 
\citet[][\S8.2 and references therein]{MR595165}, 
there are several approaches of 
establishing asymptotic normality
of $\widehat{\bm\mu}$ 
depending upon the various 
scenarios of the weights generating
function $J$ and the underlying cdf $F$.

In particular,
\citet[][{\sc{Remark}} 9]{MR0203874}
established that the $k$-variate vector
$\sqrt{n}(\widehat{\bm\mu}-\bm{\mu})$,
converges in distribution to the 
$k$-variate normal random vector
with mean $\mathbf{0}$ and 
the variance-covariance matrix $\mathbf{\Sigma}:= \left[\sigma_{ij}^{2}\right]_{i,j=1}^{k}$ 
with the entries 
\begin{align}
\sigma_{ij}^{2} 
& = 
\int_{0}^{1} 
{\alpha_{i}(u)\alpha_{j}(u) \, du}  
\label{eqn:mtm_var_cov1} \\[5pt]
& = 
\int_{0}^{1}
\left(\dfrac{1}{\ol{u}} \right)^{2}
\left[\int_{u}^{1}{\ol{v}g_{i}(v)} \, dv 
\int_{u}^{1}{\ol{w}g_{j}(w)} \, dw\right] du
\label{eqn:mtm_var_cov2} \\[5pt]
& = 
\int_{0}^{1} 
\int_{0}^{1} 
g_{i}(v) 
g_{j}(w)
K(v,w) \, dv \, dw,
\label{eqn:mtm_var_cov3}
\end{align}
where the function $K(v,w)$ is defined as
\begin{align*}
K(v,w) 
& := 
K(w, v)
=
v \wedge w 
- 
vw
= 
\min \{v, w \} - vw,
\quad 
\mbox{for} 
\quad 
0 \le v, w \le 1.
\end{align*}

By using "Integration by Parts" 
and for $0 \le a < b \le 1$,
we note that 
\begin{align}
\label{eqn:IntOfH2}
I(a,b) 
& : = 
\int_{a}^{b}
v H'(v) \, dv 
=
\left. 
v H(v)
\right|_{v = a}^{v = b}
-
\int_{a}^{b} H(v) \, dv
=
b H(b) 
-
a H(a)
- 
\int_{a}^{b} H(v) \, dv, \\[5pt]
\ol{I}(a,b)
\label{eqn:IntOfH1}
& : = 
\int_{a}^{b}
\ol{v} H'(v) \, dv 
= 
\left. 
\ol{v} \ H(v)
\right|_{v = a}^{v = b}
+
\int_{a}^{b} H(v) \, dv
=
\ol{b} H(b) 
-
\ol{a} H(a)
+ 
\int_{a}^{b} H(v) \, dv. 
\end{align} 

Further, it clearly follows that 
\begin{align*}
\int_{a}^{b}
\left( 
\dfrac{1}{\ol{u}}
\right)^{2} \, du
& = 
\left. 
\dfrac{1}{\ol{u}}
\right|_{a}^{b}
= 
\dfrac{1}{\ol{b}}
-
\dfrac{1}{\ol{a}}
=
\dfrac{\ol{a} - \ol{b}}{\ol{a} \ \ol{b}}
=
\dfrac{b - a}{\ol{a} \ \ol{b}}.
\end{align*}

\label{Page:Motivation1}
The primary motivation of this scholarly paper is
to establish computational formulas for the 
asymptotic variance,
given by Eq. \eqref{eqn:mtm_var_cov3},
of two robust $L$-estimators: 
the {\em method of trimmed moments} (MTM) and 
the {\em method of winsorized moments} (MWM). 
We also demonstrate that two asymptotic approaches
for MTM are equivalent for a specific choice of
the weights-generating function.
Our goal is to enhance the versatility of these estimators, 
allowing their application across diverse underlying distributions, 
including specific actuarial contexts such as 
payment-per-payment and payment-per-loss data scenarios 
\citep[see, e.g.,][Section 8.2 on p. 126]{MR3890025}. 
The effectiveness of our approach depends on the
availability of the cumulative distribution function, 
which broadens the applicability of our formulas and
allows us to achieve the desired level of robustness 
through data trimming or winsorization.

\section{Method of Trimmed Moments}
\label{sec:MTM}

As before, 
let 
$X_{1:n} \leq X_{2:n} \leq \cdots \leq X_{n:n}$
be the order statistics of a random sample 
$X_{1}, X_{2},..., X_{n} \stackrel{i.i.d.}{\sim} X$,
where $X \sim F(x|\bm{\theta})$ 
with $k$ unknown parameters 
$\bm{\theta} = (\theta_{1},\dots,\theta_{k})$.
Then the {\em method of trimmed moments} (MTM)
estimators of 
$\theta_{1}, \theta_{2},...,\theta_{k}$
are found as follows:

\begin{itemize}
\item 
Compute the sample trimmed moments 
\begin{equation} 
\label{eq:sample_mtmG}
\widehat{\mu}_{j} 
=
\frac{1}{n-\lfloor na_{j} \rfloor - \lfloor nb_{j} \rfloor}
\sum_{i=\lfloor na_{j} \rfloor+1}^{n-\lfloor nb_{j} \rfloor}{h_{j}(X_{i:n})},
\quad 1 \leq j \leq k.
\end{equation}

The $h_{j}'s$ in (\ref{eq:sample_mtmG}) 
are specially chosen functions, 
where the proportions 
$0 \le a_{j}, b_{j} \le 1$ 
are chosen by researcher.

\item
Compute the corresponding population trimmed moments 
\begin{equation} 
\label{eq:pop_mtmG}
\mu_{j} 
=
\frac{1}{1-a_{j}-b_{j}}\int_{a_{j}}^{\ol{b}_{j}}{h_{j}
\lp 
F^{-1}(u|\bm{\theta})
\rp \, du,\ \ \ 1\leq{j}\leq{k}}.
\end{equation} 
In (\ref{eq:pop_mtmG}), $F^{-1}(u|\bm{\theta})=\inf\:\{x:F(x|\bm{\theta})\geq{u}\}$ is the quantile function.

\item 
Now, match the sample and population trimmed moments from (\ref{eq:sample_mtmG}) and (\ref{eq:pop_mtmG}) to get the following system of equations for $\theta_{1},\theta_{2},...,\theta_{k}$
\begin{equation} 
\label{eq:match_mtm} 
\left\{
\begin{array}{lcl}
\mu_1 (\theta_1, \ldots, \theta_k) 
& = & 
\widehat{\mu}_1 \\
& \vdots & \\
\mu_k (\theta_1, \ldots, \theta_k) 
& = & 
\widehat{\mu}_k  \\
\end{array} \right.
\end{equation}
\end{itemize}

A solution, say 
$\widehat{\bm{\theta}}
=
\lp 
\widehat{\theta}_{1},\widehat{\theta}_{2},...,\widehat{\theta}_{k}
\rp$, 
if it exists, 
to the system of equations (\ref{eq:match_mtm}) 
is called the {\textit{method of trimmed moments (MTM)}} 
estimator of $\bm{\theta}$. 
Thus, $\widehat{\theta}_{j}
=:
\tau_{j}
\left(
\widehat{\mu}_{1},\widehat{\mu}_{2},..., \widehat{\mu}_{k}
\right)$, $1\leq{j}\leq{k}$, 
are the MTM estimators of 
$\theta_{1},\theta_{2},...,\theta_{k}$.

Asymptotically, Eq. \eqref{eq:sample_mtmG},
\label{Page:Weight1}
see \citet[][p. 264]{MR595165},
with the condition 
$0 \le a_{j} < \ol{b}_{j} \le 1$ 
such that $a_{j} + b_{j} < 1$, 
is equivalent to Eq. \eqref{eqn:S1} 
with the specified weights-generating function:
\begin{align}
\label{eqn:J_Fun1}
J_{j}(s) 
& =
\left\{ 
\begin{array}{ll}
(1-a_{j}-b_{j})^{-1}; & a_{j} < s < \ol{b}_{j}, \\
0; & \text{otherwise}. \\
\end{array}
\right. 
\quad 
1 \le j \le k,
\end{align}
Clearly, 
Eq. \eqref{eq:pop_mtmG}
is also equivalent to
Eq. \eqref{eqn:P1}.

With the weights-generating
function given by Eq. \eqref{eqn:J_Fun1}, 
\cite{MR2416885,MR2497558} independently 
established an equivalent simplified 
expression for $\sigma_{ij}^{2}$, 
which can be written as
\begin{align} 
\label{eqn:mtmVB_var_cov1}
\sigma_{ij}^{2}
& =
\Gamma 
\int_{a_{i}}^{\ol{b}_{i}}\int_{a_{j}}^{\ol{b}_{j}}
{K(v,w)} \ dh_{j}(F^{-1}(v)) \ dh_{i}(F^{-1}(w))
=
\Gamma 
\int_{a_{i}}^{\ol{b}_{i}}\int_{a_{j}}^{\ol{b}_{j}}
{K(v,w)} 
H_{j}'(v) H_{i}'(w) \, dv \, dw,
\end{align}
where 
\(
\displaystyle
\Gamma 
: = 
\prod_{r = i, j}
\left(
1-a_{r}-b_{r}
\right)^{-1}.
\)

Before presenting the main result of this section, 
we first state Fubini's Theorem without delving 
into the concept of product measure.

\begin{thm}[Fubini's Theorem - \citet{MR1681462}, p. 67]
\label{thm:Fubini1}
Let $E := [a,b] \times [c,d] \subseteq \R^{2}$,
s.t. $a < b$ and $c < d$, 
is a rectangular region.
Consider a Lebesgue measurable function
\[
f : E \to \R
\quad 
\mbox{s.t.}
\quad 
\int_{E} \left| f(x,y) \right| \, dy \, dx 
< 
\infty,
\]
then the Fubini's Theorem states that
\begin{align*}
\int_{a}^{b}
\int_{c}^{d} f(x,y) \, dy \, dx
& = 
\int_{E_{x}}
\left( 
\int_{E_{y}}
f(x,y) \, dy 
\right) \, dx 
=
\int_{E_{y}}
\left( 
\int_{E_{x}}
f(x,y) \, dx 
\right) \, dy, 
\\[5pt]
\mbox{i.e.,} \quad 
\int_{E} f(x,y) \, dy \, dx
& = 
\int_{E_{x}}
\left( 
\int_{E_{y}}
f_{x}(y) \, dy 
\right) \, dx 
=
\int_{E_{y}}
\left( 
\int_{E_{x}}
f_{y}(x) \, dx 
\right) \, dy,
\end{align*}
where \citep[][p. 65]{MR1681462}
\begin{align*}
\mbox{\bf $x$-section}: 
& \quad  
E_{x} 
:= 
\left\{ 
y \in \R : (x,y) \in E
\right\} 
\quad 
\mbox{and}
\quad 
f_{x}(y) 
:= 
f(x,y),
\\ 
\mbox{\bf $y$-section}: 
& \quad  
E_{y} 
:= 
\left\{ 
x \in \R : (x,y) \in E
\right\}
\quad 
\mbox{and}
\quad  
f_{y}(x) 
:= 
f(x,y).
\end{align*}
\end{thm}

\label{Section2Eliminated}

An immediate application of the 
Fubini's Theorem is the following corollary. 

\begin{cor}
\label{prop:sDoubleInt1}
Let $g_{1} : \R \to \R$ and $g_{2} : \R \to \R$ 
are two Lebesgue integrable functions 
in the finite interval $[a, b]$
with $a < b$ such that 
\[
\int_{E} 
\left| g_{1}(x) g_{2}(y) \right|
\, dy \, dx 
< 
\infty,
\]
where $E = [a, b] \times [a, b]$.
Then it follows that 
\begin{align}
\label{eqn:sDoubleInt2}
\left( 
\int_{a}^{b}
g_{1}(x) \, dx 
\right)
\left( 
\int_{a}^{b}
g_{2}(y) \, dy 
\right)
& = 
\int_{E}
g_{1}(x) g_{2}(y) \, dy \, dx
=
\int_{a}^{b}
\int_{a}^{b}
g_{1}(x) g_{2}(y) \, dy \, dx.
\end{align}
In particular if 
$g_{1} \equiv g_{2} = g$,
then Eq. \eqref{eqn:sDoubleInt2}
takes the form 
\begin{align*}
\left( 
\int_{a}^{b}
g(x) \, dx 
\right)^{2} 
& =
\left( 
\int_{a}^{b}
g(x) \, dx 
\right)
\left( 
\int_{a}^{b}
g(y) \, dy 
\right)
= 
\int_{a}^{b}
\int_{a}^{b}
g(x) g(y) \, dy \, dx.
\end{align*}
\end{cor}

With Theorem \ref{thm:Fubini1} and 
Corollary \ref{prop:sDoubleInt1},
we are now ready to establish one 
of the main results of this scholarly work.

\begin{thm}
\label{thm:MTM_Theorem1}
Equations 
\eqref{eqn:mtm_var_cov1},
\eqref{eqn:mtm_var_cov2},
\eqref{eqn:mtm_var_cov3}, and
\eqref{eqn:mtmVB_var_cov1}
are equivalent.

\begin{proof} 
Eq. \eqref{eqn:mtm_var_cov2}
follows obviously from 
\eqref{eqn:mtm_var_cov1}.

Eq. \eqref{eqn:mtm_var_cov2}
implies 
Eq. \eqref{eqn:mtm_var_cov3}:
From Corollary \ref{prop:sDoubleInt1},
we have
\begin{align}
\alpha_{i}(u)\alpha_{j}(u)
& =
\left(\dfrac{1}{\ol{u}} \right)^{2}
\left[
\int_{u}^{1}{
\ol{v}g_{i}(v)} \, dv 
\int_{u}^{1}{\ol{w}g_{j}(w)} \, dw
\right] 
\nonumber \\[5pt]
& = 
\left(\dfrac{1}{\ol{u}} \right)^{2}
\int_{u}^{1}
\int_{u}^{1}
\ol{v} \ \ol{w} \ 
g_{i}(v) g_{j}(w) \, dv \, dw.
\nonumber
\end{align}

\begin{figure}[hbt!]
\centering
\includegraphics[width=0.85\textwidth] 
{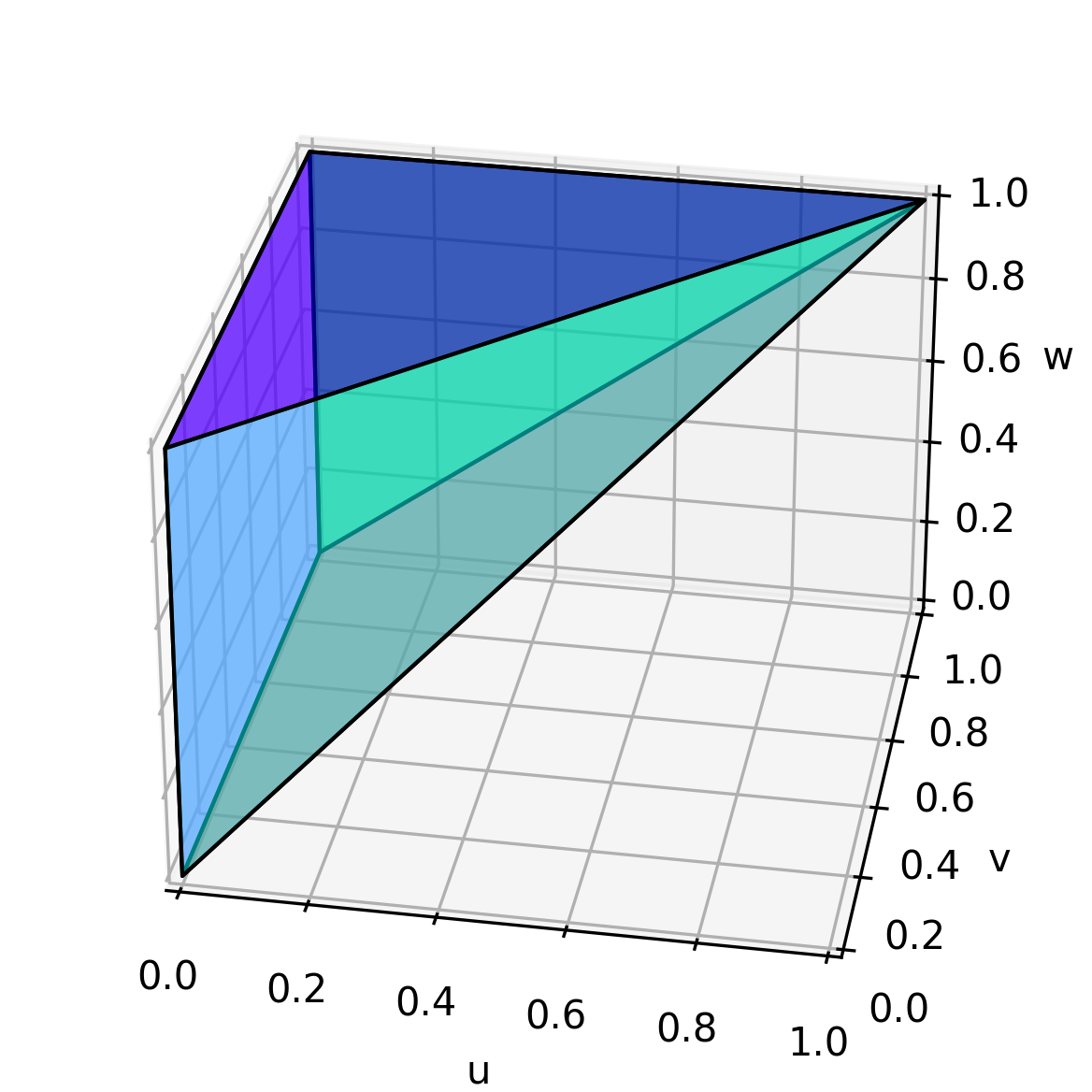}

\vspace{-0.25cm}

\caption{Visualization of the rectangular pyramid $\Omega$.}
\label{fig:figure1}
\end{figure}

Then from Eq. \eqref{eqn:mtm_var_cov2},
it follows that 
\begin{align}
\sigma_{ij}^{2} 
& = 
\int_{0}^{1}
\left(\dfrac{1}{\ol{u}} \right)^{2}
\left[ 
\int_{u}^{1}
\int_{u}^{1}
\ol{v} \ \ol{w} \ 
g_{i}(v) g_{j}(w) \, dv \, dw
\right] \, du 
\nonumber \\
& = 
\int_{\Omega}
\left(\dfrac{1}{\ol{u}} \right)^{2}
\ol{v} \ \ol{w} \ 
g_{i}(v) g_{j}(w)
\, dv \, dw \, du,
\label{eqn:mtm_var_cov4} 
\end{align}
where the domain of integration, 
$\Omega \subset [0,1]^{3}$, 
in Eq. \eqref{eqn:mtm_var_cov4}, 
is defined as
\begin{align*}
\Omega 
& := 
\left\{
0 \le u \le 1, \ 
u \le v \le 1, \ 
u \le w \le 1
\right\} 
\nonumber \\ 
& =
\left\{ 
0 \le x \le v \wedge w, \ 
0 \le v \le 1, \ 
0 \le w \le 1, 
\right\}
\subseteq \R^{3}.
\end{align*}

It is clear that 
$\Omega$
--
a rectangular pyramid as shown 
in Figure \ref{fig:figure1},
is a measurable subset of the unit
cube $[0,1]^{3}$. 
That is,
the Lebesgue measure of 
$\Omega$ is simply the
volume of $\Omega$ which is 1/3.
Also, note that 
$\Omega_{u} = v \wedge w$, 
$\Omega_{v} = [0,1]$, and 
$\Omega_{w} = [0,1]$.
Since $\Omega$ is measurable 
and assuming the functions 
involved in Eq. \eqref{eqn:mtm_var_cov4} 
satisfy the conditions for Fubini's Theorem,
i.e., Corollary \ref{prop:sDoubleInt1},
then the order of integration in 
Eq. \eqref{eqn:mtm_var_cov4} can be changed.
That is, 

\begin{align}
\sigma_{ij}^{2} 
& =
\int_{\Omega}
\left(\dfrac{1}{\ol{u}} \right)^{2}
\ol{v} \ \ol{w} \ 
g_{i}(v) g_{j}(w)
\, dv \, dw \, du 
\nonumber \\ 
& = 
\int_{\Omega_{u}}
\int_{\Omega_{w}}
\int_{\Omega_{v}}
\left(\dfrac{1}{\ol{u}} \right)^{2}
\ol{v} \ \ol{w} \ 
g_{i}(v) g_{j}(w)
\, dv \, dw \, du 
\nonumber \\ 
& = 
\int_{\Omega_{w}}
\int_{\Omega_{v}}
\int_{\Omega_{u}}
\left(\dfrac{1}{\ol{u}} \right)^{2}
\ol{v} \ \ol{w} \ 
g_{i}(v) g_{j}(w)
\, du \, dv \, dw 
\nonumber \\ 
& = 
\int_{\Omega_{w}}
\int_{\Omega_{v}}
\ol{v} \ \ol{w} \ 
g_{i}(v) g_{j}(w)
\left( 
\int_{\Omega_{u}}
\left(\dfrac{1}{\ol{u}} \right)^{2} \, du
\right) \, dv \, dw 
\nonumber \\ 
& = 
\int_{\Omega_{w}}
\int_{\Omega_{v}}
\ol{v} \ \ol{w} \ 
g_{i}(v) g_{j}(w)
\left. 
\left( 
\frac{1}{\ol{u}}
\right)
\right|_{u=0}^{u = v \wedge w} 
\, dv \, dw 
\nonumber \\ 
& = 
\int_{\Omega_{w}}
\int_{\Omega_{v}}
\ol{v} \ \ol{w} \ 
g_{i}(v) g_{j}(w)
\left( 
\frac{1}{1 - v \wedge w}
-
1
\right) \, dv \, dw. 
\label{eqn:mtm_var_cov5}
\end{align}
As observed in Eq. \eqref{eqn:mtm_var_cov5},
in order to establish 
Eq. \eqref{eqn:mtm_var_cov2}
implies 
Eq. \eqref{eqn:mtm_var_cov3}, 
it simply remains to show that 
\begin{align*}
\label{eqn:kFun3}
K(v,w) 
& = 
v \wedge w 
- 
vw 
=
\ol{v} \ \ol{w}
\left( 
\frac{1}{1 - v \wedge w}
-
1
\right).
\end{align*}
To establish this equality, 
we investigate two scenarios.
\label{DeleteIn}
\begin{itemize}
\item[(i)] 
If $v \le w$,
then $v \wedge w = v$
and we get 
\[
\ol{v} \ \ol{w}
\left( 
\frac{1}{1 - v \wedge w}
-
1
\right) 
=
\ol{v} \ \ol{w}
\left( 
\frac{1}{1 - v}
-
1
\right) 
=
v - vw
=
K(v,w).
\]

\item[(ii)]
If $v > w$,
then $v \wedge w = v$
and it follows that
\[
\ol{v} \ \ol{w}
\left( 
\frac{1}{1 - v \wedge w}
-
1
\right) 
=
\ol{v} \ \ol{w}
\left( 
\frac{1}{1 - w}
-
1
\right) 
=
w - vw
=
K(v,w).
\]
\end{itemize}

Therefore, 
Eq. \eqref{eqn:mtm_var_cov5}
takes the form
\begin{align}
\sigma_{ij}^{2} 
& =
\int_{\Omega_{w}}
\int_{\Omega_{v}}
g_{i}(v) g_{j}(w)
\left( 
\frac{1}{1 - v \wedge w}
-
1
\right) 
\ol{v} \ \ol{w} \, dv \, dw 
\nonumber \\ 
& = 
\int_{\Omega_{w}}
\int_{\Omega_{v}}
g_{i}(v) g_{j}(w)
K(v,w) \, dv \, dw 
\nonumber \\ 
& = 
\int_{0}^{1}
\int_{0}^{1}
g_{i}(v) g_{j}(w)
K(v,w) \, dv \, dw,
\nonumber
\end{align}
which is
Eq. \eqref{eqn:mtm_var_cov3}
as desired.

Eq. \eqref{eqn:mtm_var_cov3}
implies 
Eq. \eqref{eqn:mtmVB_var_cov1}:
\label{Page:Correction1}
Note the relation \eqref{eqn:cgj7}.
For the special weight function specified in Eq. \eqref{eqn:J_Fun1}, 
Eq. \eqref{eqn:mtm_var_cov3} converts to:
\begin{align}
\sigma_{ij}^{2}
& = 
\int_{0}^{1} 
\int_{0}^{1} 
g_{i}(v) 
g_{j}(w)
K(v,w) \, dv \, dw 
\nonumber \\
& = 
\int_{0}^{1} 
\int_{0}^{1} 
J_{j}(v) J_{i}(w) 
K(v,w) \, 
H_{j}'(v) H_{i}'(w) \, dv \, dw
\nonumber \\
& = 
\Gamma
\int_{a_{i}}^{\ol{b}_{i}} 
\int_{a_{j}}^{\ol{b}_{j}} 
K(v,w) \, 
H_{j}'(v) H_{i}'(w) \, dv \, dw  
\nonumber \\ 
& = 
\Gamma
\int_{a_{i}}^{\ol{b}_{i}} 
\int_{a_{j}}^{\ol{b}_{j}} 
K(v,w)
dh_{j}(F^{-1}(v)) \, 
dh_{i}(F^{-1}(w)),
\label{eqn:cgj9}
\end{align}
which is exactly Eq. \eqref{eqn:mtmVB_var_cov1}.
\end{proof}
\end{thm}

We now evaluate the expression 
for $\sigma_{ij}^{2}$ from 
Eq. \eqref{eqn:cgj9} in more explicit form. 
Altogether there are six 
possible combinations of 
trimming proportions 
$(a_{i},b_{i})$ and $(a_{j},b_{j})$
for $1 \le i, j \le k$, 
and they are: 
\begin{enumerate}[label=(\roman*)]
\item 
$0 \leq a_{i} \leq a_{j} \leq \ol{b}_{i} \leq \ol{b}_{j} \leq 1$,
{\scriptsize
\setlength{\unitlength}{1cm}
\thinlines
\begin{picture}(5,0.2)
\put(4,.2){\vector(1,0){3.5}}
\put(4,.2){\vector(-1,0){3}}
\put(1.5,.2){\line(0,1){.2}} 
\put(2.5,.2){\line(0,1){.2}}
\put(3.5,.2){\line(0,1){.2}}
\put(4.6,.2){\line(0,1){.2}}
\put(6.1,.2){\line(0,1){.2}}
\put(7.0,.2){\line(0,1){.2}}
\put(1.43,-.1){$0$}
\put(2.40,-.1){$a_{i}$}
\put(3.40,-.1){$a_{j}$}
\put(4.50,-.1){$\ol{b}_{i}$}
\put(6.05,-.1){$\ol{b}_{j}$}
\put(6.93,-.1){$1$}
\end{picture} 
}

\item 
$0 \leq a_{i} \leq \ol{b}_{i} \leq a_{j} \leq \ol{b}_{j} \leq 1$,
{\scriptsize
\setlength{\unitlength}{1cm}
\thinlines
\begin{picture}(5,0.2)
\put(4,.2){\vector(1,0){3.5}}
\put(4,.2){\vector(-1,0){3}}
\put(1.5,.2){\line(0,1){.2}}
\put(2.5,.2){\line(0,1){.2}}
\put(3.5,.2){\line(0,1){.2}}
\put(4.6,.2){\line(0,1){.2}}
\put(6.1,.2){\line(0,1){.2}}
\put(7.0,.2){\line(0,1){.2}}
\put(1.43,-.1){$0$}
\put(2.40,-.1){$a_{i}$}
\put(3.40,-.1){$\ol{b}_{i}$}
\put(4.50,-.1){$a_{j}$}
\put(6.05,-.1){$\ol{b}_{j}$}
\put(6.93,-.1){$1$}
\end{picture} 
}

\item $0 \leq a_{i} \leq a_{j} \leq \ol{b}_{j} \leq \ol{b}_{i} \leq 1$,
{\scriptsize
\setlength{\unitlength}{1cm}
\thinlines
\begin{picture}(5,0.2)
\put(4,.2){\vector(1,0){3.5}}
\put(4,.2){\vector(-1,0){3}}
\put(1.5,.2){\line(0,1){.2}} 
\put(2.5,.2){\line(0,1){.2}}
\put(3.5,.2){\line(0,1){.2}}
\put(4.6,.2){\line(0,1){.2}}
\put(6.1,.2){\line(0,1){.2}}
\put(7.0,.2){\line(0,1){.2}}
\put(1.43,-.1){$0$}
\put(2.40,-.1){$a_{i}$}
\put(3.40,-.1){$a_{j}$}
\put(4.50,-.1){$\ol{b}_{j}$}
\put(6.05,-.1){$\ol{b}_{i}$}
\put(6.93,-.1){$1$}
\end{picture} 
}

\item $0 \leq a_{j} \leq \ol{b}_{j} \leq a_{i} \leq \ol{b}_{i} \leq 1$,
{\scriptsize
\setlength{\unitlength}{1cm}
\thinlines
\begin{picture}(5,0.2)
\put(4,.2){\vector(1,0){3.5}}
\put(4,.2){\vector(-1,0){3}}
\put(1.5,.2){\line(0,1){.2}}
\put(2.5,.2){\line(0,1){.2}}
\put(3.5,.2){\line(0,1){.2}}
\put(4.6,.2){\line(0,1){.2}}
\put(6.1,.2){\line(0,1){.2}}
\put(7.0,.2){\line(0,1){.2}}
\put(1.43,-.1){$0$}
\put(2.40,-.1){$a_{j}$}
\put(3.40,-.1){$\ol{b}_{j}$}
\put(4.50,-.1){$a_{i}$}
\put(6.05,-.1){$\ol{b}_{i}$}
\put(6.93,-.1){$1$}
\end{picture} 
}

\item $0 \leq a_{j} \leq a_{i} \leq \ol{b}_{j} \leq \ol{b}_{i} \leq 1$,
{\scriptsize
\setlength{\unitlength}{1cm}
\thinlines
\begin{picture}(5,0.2)
\put(4,.2){\vector(1,0){3.5}}
\put(4,.2){\vector(-1,0){3}}
\put(1.5,.2){\line(0,1){.2}} 
\put(2.5,.2){\line(0,1){.2}}
\put(3.5,.2){\line(0,1){.2}}
\put(4.6,.2){\line(0,1){.2}}
\put(6.1,.2){\line(0,1){.2}}
\put(7.0,.2){\line(0,1){.2}}
\put(1.43,-.1){$0$}
\put(2.40,-.1){$a_{j}$}
\put(3.40,-.1){$a_{i}$}
\put(4.50,-.1){$\ol{b}_{j}$}
\put(6.05,-.1){$\ol{b}_{i}$}
\put(6.93,-.1){$1$}
\end{picture} 
}
\item $0 \leq a_{j} \leq a_{i} \leq \ol{b}_{i} \leq \ol{b}_{j} \leq 1$.
{\scriptsize
\setlength{\unitlength}{1cm}
\thinlines
\begin{picture}(5,0.2)
\put(4,.2){\vector(1,0){3.5}}
\put(4,.2){\vector(-1,0){3}}
\put(1.5,.2){\line(0,1){.2}} 
\put(2.5,.2){\line(0,1){.2}}
\put(3.5,.2){\line(0,1){.2}}
\put(4.6,.2){\line(0,1){.2}}
\put(6.1,.2){\line(0,1){.2}}
\put(7.0,.2){\line(0,1){.2}}
\put(1.43,-.1){$0$}
\put(2.40,-.1){$a_{j}$}
\put(3.40,-.1){$a_{i}$}
\put(4.50,-.1){$\ol{b}_{i}$}
\put(6.05,-.1){$\ol{b}_{j}$}
\put(6.93,-.1){$1$}
\end{picture} 
}
\end{enumerate}

Each choice results in a robust estimator
with different robustness and efficiency properties. 
Without loss of generality, 
we proceed with scenario (i) only
in this manuscript. 
The other scenarios can be handled similarly.
That is, we evaluate $\sigma_{ij}^{2}$ 
given by Eq. \eqref{eqn:mtmVB_var_cov1}
when the trimming proportions 
satisfy the following inequality:
\begin{align}
\label{eqn:abCondition1}
0 
\le 
a_{i} 
\le 
a_{j} 
< 
\ol{b}_{i} 
\le 
\ol{b}_{j} 
\le 
1.
\end{align}

\begin{thm}
\label{thm:MTM_Var1}
With the trimming proportions satisfying 
inequality \eqref{eqn:abCondition1}, 
the variance-covariance entry 
$\sigma_{ij}^{2}$ of Eq. \eqref{eqn:mtmVB_var_cov1}
is given by 
\(
\sigma_{ij}^{2}
= 
\Gamma V_{11},
\)
where 
\begin{eqnarray*}
V_{11}
& := &
I_{i}(a_{i},a_{j}) \, 
\ol{I}_{j}(a_{j},\ol{b}_{j})
+
b_{j} H_{j}(\ol{b}_{j}) 
I_{i}(a_{j},\ol{b}_{i})
-
a_{j} H_{j}(a_{j})
\ol{I}_{i}(a_{j},\ol{b}_{i})
\nonumber \\
& & 
- 
b_{i} H_{i}(\ol{b}_{i})
\int_{a_{j}}^{\ol{b}_{i}}
H_{j}(v) \, dv 
+ 
\int_{a_{j}}^{\ol{b}_{i}}
H_{i}(v) H_{j}(v) \, dv 
\nonumber \\ 
& & 
- 
a_{j} H_{i}(a_{j})
\int_{a_{j}}^{\ol{b}_{i}} 
H_{j}(v) \, dv 
-
\int_{a_{j}}^{\ol{b}_{i}}
\int_{v}^{\ol{b}_{i}}
H_{i}(w)
H_{j}(v) \, dw \, dv
-
\int_{a_{j}}^{\ol{b}_{i}} 
\int_{a_{j}}^{v}
H_{i}(w)
H_{j}(v) \, dw \, dv
\nonumber \\
& & 
+ 
\ol{b}_{i}H_{i}(\ol{b}_{i})
\int_{\ol{b}_{i}}^{\ol{b}_{j}}
H_{j}(v) \, dv 
-
a_{j} 
H_{j}(a_{j})
\int_{\ol{b}_{i}}^{\ol{b}_{j}}
H_{j}(v) \, dv
-
\int_{a_{j}}^{\ol{b}_{i}}
\int_{\ol{b}_{i}}^{\ol{b}_{j}}
H_{i}(w) H_{j}(v) \, dv \, dw.
\end{eqnarray*}

\begin{proof}
From Eq. \eqref{eqn:mtmVB_var_cov1}, 
we have
\begin{align*} 
\sigma_{ij}^{2}
& =
\Gamma 
\int_{a_{i}}^{\ol{b}_{i}} 
\int_{a_{j}}^{\ol{b}_{j}} 
K(v,w)
H_{j}^{'}(v)
H_{i}^{'}(w) \, dv \, dw.
\end{align*}

Thus, we just need to evaluate
\begin{eqnarray}
V_{11}
& = & 
\int_{a_{i}}^{\ol{b}_{i}} 
\int_{a_{j}}^{\ol{b}_{j}} 
K(v,w)
H_{j}^{'}(v)
H_{i}^{'}(w) \, dv \, dw
\nonumber \\ 
& = & 
\int_{a_{i}}^{\ol{b}_{i}} 
\int_{a_{j}}^{\ol{b}_{j}} 
\left[ 
v \wedge w - vw
\right]
H_{j}^{'}(v)
H_{i}^{'}(w) \, dv \, dw.
\nonumber \\
& = & 
\int_{a_{i}}^{a_{j}} 
\int_{a_{j}}^{\ol{b}_{j}} 
\left[ 
v \wedge w - vw
\right]
H_{j}^{'}(v)
H_{i}^{'}(w) \, dv \, dw
+ 
\int_{a_{j}}^{\ol{b}_{i}} 
\int_{a_{j}}^{\ol{b}_{j}} 
\left[ 
v \wedge w - vw
\right]
H_{j}^{'}(v)
H_{i}^{'}(w) \, dv \, dw
\nonumber \\ 
& = & 
\int_{a_{i}}^{a_{j}} 
\int_{a_{j}}^{\ol{b}_{j}} 
w \, \ol{v}
H_{j}^{'}(v)
H_{i}^{'}(w) \, dv \, dw
+ 
\int_{a_{j}}^{\ol{b}_{i}} 
\int_{a_{j}}^{\ol{b}_{j}} 
\left[ 
v \wedge w - vw
\right]
H_{j}^{'}(v)
H_{i}^{'}(w) \, dv \, dw
\nonumber \\ 
& = & 
\lp 
\int_{a_{i}}^{a_{j}} 
w H_{i}^{'}(w) \, dw 
\rp 
\lp 
\int_{a_{j}}^{\ol{b}_{j}} 
\ol{v}
H_{j}^{'}(v) \, dv 
\rp 
+ 
\int_{a_{j}}^{\ol{b}_{i}} 
\int_{a_{j}}^{\ol{b}_{j}} 
\left[ 
v \wedge w - vw
\right]
H_{j}^{'}(v)
H_{i}^{'}(w) \, dv \, dw
\nonumber \\
& = & 
I_{i}(a_{i},a_{j}) \, 
\ol{I}_{j}(a_{j},\ol{b}_{j}) 
+ 
V_{11}^{*}
\label{eqn:cgj11} \\ 
V_{11}^{*}
& := & 
\int_{a_{j}}^{\ol{b}_{i}} 
\int_{a_{j}}^{\ol{b}_{j}} 
\left[ 
v \wedge w - vw
\right]
H_{j}^{'}(v)
H_{i}^{'}(w) \, dv \, dw 
\nonumber \\ 
& = & 
\int_{a_{j}}^{\ol{b}_{i}} 
H_{i}^{'}(w) 
\left[ 
\int_{a_{j}}^{w} 
\left[ 
v \wedge w - vw
\right]
H_{j}^{'}(v) \, dv
+ 
\int_{w}^{\ol{b}_{j}} 
\left[ 
v \wedge w - vw
\right]
H_{j}^{'}(v) \, dv
\right] dw 
\nonumber \\ 
& = & 
\int_{a_{j}}^{\ol{b}_{i}} 
H_{i}^{'}(w) 
\left[ 
\int_{a_{j}}^{w} 
\left[ 
v - vw
\right]
H_{j}^{'}(v) \, dv
+ 
\int_{w}^{\ol{b}_{j}} 
\left[ 
w - vw
\right]
H_{j}^{'}(v) \, dv
\right] dw
\nonumber \\ 
& = & 
\int_{a_{j}}^{\ol{b}_{i}} 
\ol{w}
H_{i}^{'}(w) 
\left[ 
\int_{a_{j}}^{w} 
v H_{j}^{'}(v) \, dv 
\right] dw 
+ 
\int_{a_{j}}^{\ol{b}_{i}} 
w
H_{i}^{'}(w) 
\left[ 
\int_{w}^{\ol{b}_{j}} 
\ol{v}
H_{j}^{'}(v) \, dv
\right] dw 
\nonumber \\ 
& = & 
\int_{a_{j}}^{\ol{b}_{i}} 
\ol{w}
H_{i}^{'}(w) 
\lb 
w H_{j}(w) 
- 
a_{j} H_{j}(a_{j}) 
-
\int_{a_{j}}^{w} H_{j}(v) \, dv
\rb dw 
\quad 
\lp \mbox{from Eq. \eqref{eqn:IntOfH1}} \rp
\nonumber \\ 
& & 
+ 
\int_{a_{j}}^{\ol{b}_{i}} 
w
H_{i}^{'}(w) 
\lb 
b_{j} H_{j}(\ol{b}_{j}) 
- 
\ol{w}  H_{j}(w) 
+ 
\int_{w}^{\ol{b}_{j}} H_{j}(v) \, dv
\rb dw
\quad 
\lp \mbox{from Eq. \eqref{eqn:IntOfH2}} \rp
\nonumber \\ 
& = & 
-
a_{j} H_{j}(a_{j})
\int_{a_{j}}^{\ol{b}_{i}} 
\ol{w} 
H_{i}^{'}(w) \, dw
-
\int_{a_{j}}^{\ol{b}_{i}} 
\int_{a_{j}}^{w}
\ol{w}
H_{j}(v) H_{i}^{'}(w) \, dv \, dw 
\nonumber \\ 
& & 
+ 
b_{j} H_{j}(\ol{b}_{j}) 
\int_{a_{j}}^{\ol{b}_{i}} 
w
H_{i}^{'}(w) \, dw 
+
\int_{a_{j}}^{\ol{b}_{i}}
\int_{w}^{\ol{b}_{j}}
w
H_{j}(v)
H_{i}^{'}(w) \, dv \, dw 
\nonumber \\ 
& = & 
-
a_{j} H_{j}(a_{j})
\ol{I}_{i}(a_{j},\ol{b}_{i})
-
\int_{a_{j}}^{\ol{b}_{i}}
H_{j}(v) 
\lb 
\int_{v}^{\ol{b}_{i}}
\ol{w} H_{i}^{'}(w) \, dw
\rb dv 
+ 
b_{j} H_{j}(\ol{b}_{j}) 
I_{i}(a_{j},\ol{b}_{i}) 
\nonumber \\ 
& & 
+ 
\int_{a_{j}}^{\ol{b}_{i}} 
H_{j}(v) 
\lb 
\int_{a_{j}}^{v} 
w H_{i}^{'}(w) \, dw
\rb dv 
+ 
\underbrace{
\int_{\ol{b}_{i}}^{\ol{b}_{j}}
H_{j}(v) 
\lb 
\int_{a_{j}}^{\ol{b}_{i}}
w H_{i}^{'}(w) \, dw
\rb dv
}_{=:A_{1}}
\label{eqn:cgj12} \\ 
& = & 
-
a_{j} H_{j}(a_{j})
\ol{I}_{i}(a_{j},\ol{b}_{i})
+
b_{j} H_{j}(\ol{b}_{j}) 
I_{i}(a_{j},\ol{b}_{i})
\nonumber \\ 
& & 
-
\int_{a_{j}}^{\ol{b}_{i}}
H_{j}(v)
\lb 
b_{i} H_{i}(\ol{b}_{i}) 
- 
\ol{v} H_{i}(v) 
+ 
\int_{v}^{\ol{b}_{i}} H_{i}(w) \, dw
\rb dv
\nonumber \\ 
& & 
+ 
\int_{a_{j}}^{\ol{b}_{i}} 
H_{j}(v) 
\lb 
v H_{i}(v) 
-
a_{j} H_{i}(a_{j}) 
- 
\int_{a_{j}}^{v} H_{i}(w) \, dw
\rb dv 
\nonumber \\ 
& & 
+ 
\underbrace{
\int_{\ol{b}_{i}}^{\ol{b}_{j}}
H_{j}(v) 
\lb 
\ol{b}_{i}H_{i}(\ol{b}_{i}) 
- 
a_{j} 
H_{j}(a_{j}) 
- 
\int_{a_{j}}^{\ol{b}_{i}}
H_{i}(w) \, dw
\rb dv}_{=A_{1}}
\nonumber \\ 
& = & 
-
a_{j} H_{j}(a_{j})
\ol{I}_{i}(a_{j},\ol{b}_{i})
+
b_{j} H_{j}(\ol{b}_{j}) 
I_{i}(a_{j},\ol{b}_{i}) 
\nonumber \\ 
& & 
- 
b_{i} H_{i}(\ol{b}_{i})
\int_{a_{j}}^{\ol{b}_{j}}
H_{j}(v) \, dv 
+ 
\int_{a_{j}}^{\ol{b}_{i}}
\ol{v} H_{i}(v) H_{j}(v) \, dv
-
\underbrace{
\int_{a_{j}}^{\ol{b}_{i}}
\int_{v}^{\ol{b}_{i}}
H_{i}(w)
H_{j}(v) \, dw \, dv}_{=:A_{2}} 
\nonumber \\ 
& & 
+ 
\int_{a_{j}}^{\ol{b}_{i}} 
v H_{i}(v) 
H_{j}(v) \, dv
- 
a_{j} H_{i}(a_{j})
\int_{a_{j}}^{\ol{b}_{i}} 
H_{j}(v) \, dv 
-
\underbrace{
\int_{a_{j}}^{\ol{b}_{i}} 
\int_{a_{j}}^{v}
H_{i}(w)
H_{j}(v) \, dw \, dv}_{=:A_{3}} 
\nonumber \\
& & 
+ 
\underbrace{
\int_{\ol{b}_{i}}^{\ol{b}_{j}}
H_{j}(v) 
\lb 
\ol{b}_{i}H_{i}(\ol{b}_{i}) 
- 
a_{j} 
H_{j}(a_{j}) 
- 
\int_{a_{j}}^{\ol{b}_{i}}
H_{i}(w) \, dw
\rb dv}_{=A_{1}}
\nonumber \\ 
& = & 
b_{j} H_{j}(\ol{b}_{j}) 
I_{i}(a_{j},\ol{b}_{i})
-
a_{j} H_{j}(a_{j})
\ol{I}_{i}(a_{j},\ol{b}_{i})
\nonumber \\ 
& & 
- 
b_{i} H_{i}(\ol{b}_{i})
\int_{a_{j}}^{\ol{b}_{i}}
H_{j}(v) \, dv 
+ 
\int_{a_{j}}^{\ol{b}_{i}}
H_{i}(v) H_{j}(v) \, dv 
\nonumber \\ 
& & 
- 
a_{j} H_{i}(a_{j})
\int_{a_{j}}^{\ol{b}_{i}} 
H_{j}(v) \, dv 
-
(A_{2} + A_{3}) 
\nonumber \\
& & 
+ 
\underbrace{
\ol{b}_{i}H_{i}(\ol{b}_{i})
\int_{\ol{b}_{i}}^{\ol{b}_{j}}
H_{j}(v) \, dv 
-
a_{j} 
H_{j}(a_{j})
\int_{\ol{b}_{i}}^{\ol{b}_{j}}
H_{j}(v) \, dv
-
\overbrace{
\int_{a_{j}}^{\ol{b}_{i}}
\int_{\ol{b}_{i}}^{\ol{b}_{j}}
H_{i}(w) H_{j}(v) \, dv \, dw
}^{=:A_{4}}
}_{=A_{1}}
\label{eqn:cgj13}
\end{eqnarray}

Note that 
\begin{eqnarray}
A_{2} + A_{3}
& = & 
A_{3} + A_{2}
\nonumber \\ 
& = & 
\int_{a_{j}}^{\ol{b}_{i}} 
\int_{a_{j}}^{v}
H_{i}(w)
H_{j}(v) \, dw \, dv
+ 
\int_{a_{j}}^{\ol{b}_{i}}
\int_{v}^{\ol{b}_{i}}
H_{i}(w)
H_{j}(v) \, dw \, dv
\nonumber \\ 
& = & 
\int_{a_{j}}^{\ol{b}_{i}} 
\int_{a_{j}}^{\ol{b}_{i}}
H_{i}(w)
H_{j}(v) \, dv \, dw 
=
\lp 
\int_{a_{j}}^{\ol{b}_{i}} 
H_{i}(w) \, dw 
\rp 
\lp 
\int_{a_{j}}^{\ol{b}_{i}} 
H_{j}(v) \, dv 
\rp,
\label{eqn:cgj14} \\
A_{4} 
& = & 
\int_{a_{j}}^{\ol{b}_{i}}
\int_{\ol{b}_{i}}^{\ol{b}_{j}}
H_{i}(w) H_{j}(v) \, dv \, dw 
= 
\lp 
\int_{a_{j}}^{\ol{b}_{i}}
H_{i}(w) \, dw
\rp 
\lp 
\int_{\ol{b}_{i}}^{\ol{b}_{j}}
H_{j}(v) \, dv
\rp.
\label{eqn:cgj15} 
\nonumber 
\end{eqnarray}

From Eqs. 
\eqref{eqn:cgj11},
\eqref{eqn:cgj13}, and
\eqref{eqn:cgj14},
it follows that 
\begin{eqnarray}
V_{11} 
& = & 
I_{i}(a_{i},a_{j}) \, 
\ol{I}_{j}(a_{j},\ol{b}_{j}) 
+ 
V_{11}^{*} 
\nonumber \\
& = & 
I_{i}(a_{i},a_{j}) \, 
\ol{I}_{j}(a_{j},\ol{b}_{j})
+
b_{j} H_{j}(\ol{b}_{j}) 
I_{i}(a_{j},\ol{b}_{i})
-
a_{j} H_{j}(a_{j})
\ol{I}_{i}(a_{j},\ol{b}_{i})
\nonumber \\ 
& & 
- 
b_{i} H_{i}(\ol{b}_{i})
\int_{a_{j}}^{\ol{b}_{i}}
H_{j}(v) \, dv 
+ 
\int_{a_{j}}^{\ol{b}_{i}}
H_{i}(v) H_{j}(v) \, dv 
\nonumber \\ 
& & 
- 
a_{j} H_{i}(a_{j})
\int_{a_{j}}^{\ol{b}_{i}} 
H_{j}(v) \, dv 
-
(A_{2} + A_{3}) 
\nonumber \\
& & 
+ 
\ol{b}_{i}H_{i}(\ol{b}_{i})
\int_{\ol{b}_{i}}^{\ol{b}_{j}}
H_{j}(v) \, dv 
-
a_{j} 
H_{j}(a_{j})
\int_{\ol{b}_{i}}^{\ol{b}_{j}}
H_{j}(v) \, dv
-
A_{4},
\label{eqn:cgj16}
\nonumber 
\end{eqnarray}
as desired.
\end{proof}
\end{thm}
 
For most practical estimation purposes and 
for computational simplicity, 
see, e.g., 
\cite{MR2497558, MR4263275, MR3758788}, 
we consider the following trimming proportions:
\begin{align}
\label{eqn:abSituation1}
0 
\le a = a_{i} = a_{j} 
\le
1 - b
=
1 - b_{i} = 1 - b_{j}
\le 1.
\end{align}

It is noted that \cite{MR2591318} employs trimming proportions 
\(a_{i} \neq a_{j}\) and \(b_{i} \neq b_{j}\) for \(i \neq j\); 
however, the variance-covariance formula, 
as summarized in Theorem \ref{thm:MTM_Var1},
is not presented in that paper.

\begin{cor}
\label{cor:MTM1}
For 
\(
0 \le a \le \overline{b} \le 1,
\)
define
\(
\Delta_{x} 
: = 
aH_{x}(a) 
+ 
\int_{a}^{\ol{b}} H_{x}(v) \, dv 
+ 
b H_{x}(\ol{b}), \
x \in \{ i, j \}.
\)
Then,
given the special condition specified
by inequality \eqref{eqn:abSituation1}, 
prove that 
\begin{align*}
V_{11}
& = 
a H_{i}(a) H_{j}(a) 
+
b H_{i}(\ol{b}) H_{j}(\ol{b})
+
\int_{a}^{\ol{b}}
H_{i}(v) H_{j}(v) \, dv
- 
\Delta_{i} \, \Delta_{j}.
\end{align*}

\begin{proof}
With the condition of Eq. \eqref{eqn:abSituation1}, 
from Eq. \eqref{eqn:cgj11},
it follows that 
\[
I_{i}(a_{i},a_{j}) 
= 
0
= 
\ol{I}_{j}(a_{j},\ol{b}_{j}) = 0.
\]

Similarly, from Eq. \eqref{eqn:cgj12},
it follows that 
\[
A_{1} = 0.
\]

Therefore, from Eq. \eqref{eqn:cgj16}, 
we have 
\begin{eqnarray}
V_{11} 
& = & 
b_{j} H_{j}(\ol{b}_{j}) 
I_{i}(a_{j},\ol{b}_{i})
-
a_{j} H_{j}(a_{j})
\ol{I}_{i}(a_{j},\ol{b}_{i})
\nonumber \\ 
& & 
- 
b_{i} H_{i}(\ol{b}_{i})
\int_{a_{j}}^{\ol{b}_{i}}
H_{j}(v) \, dv 
+ 
\int_{a_{j}}^{\ol{b}_{i}}
H_{i}(v) H_{j}(v) \, dv 
\nonumber \\ 
& & 
- 
a_{j} H_{i}(a_{j})
\int_{a_{j}}^{\ol{b}_{i}} 
H_{j}(v) \, dv 
-
(A_{2} + A_{3})  
\nonumber \\
& = & 
b H_{j}(\ol{b}) 
I_{i}(a,\ol{b})
-
a H_{j}(a)
\ol{I}_{i}(a,\ol{b})
\nonumber \\ 
& & 
- 
b H_{i}(\ol{b})
\int_{a}^{\ol{b}}
H_{j}(v) \, dv 
+ 
\int_{a}^{\ol{b}}
H_{i}(v) H_{j}(v) \, dv 
\nonumber \\ 
& & 
- 
a H_{i}(a)
\int_{a}^{\ol{b}} 
H_{j}(v) \, dv 
-
\lp 
\int_{a}^{\ol{b}} 
H_{i}(w) \, dw 
\rp 
\lp 
\int_{a}^{\ol{b}} 
H_{j}(v) \, dv 
\rp 
\nonumber \\ 
& = & 
b H_{j}(\ol{b}) 
\lb 
\ol{b} H_{i}(\ol{b})
- 
a H_{i}(a)
-
\int_{a}^{\ol{b}}
H_{i}(v) \, dv 
\rb
\nonumber \\ 
& & 
- 
a H_{j}(a)
\lb 
b H_{i}(\ol{b}) 
- 
\ol{a} H_{i}(a) 
+ 
\int_{a}^{\ol{b}} H_{i}(v) \, dv
\rb 
\nonumber \\
& &
- 
b H_{i}(\ol{b})
\int_{a}^{\ol{b}}
H_{j}(v) \, dv 
+ 
\int_{a}^{\ol{b}}
H_{i}(v) H_{j}(v) \, dv 
\nonumber \\ 
& & 
- 
a H_{i}(a)
\int_{a}^{\ol{b}} 
H_{j}(v) \, dv 
-
\lp 
\int_{a}^{\ol{b}} 
H_{i}(w) \, dw 
\rp 
\lp 
\int_{a}^{\ol{b}} 
H_{j}(v) \, dv 
\rp 
\nonumber \\
& = & 
a H_{i}(a) H_{j}(a) 
+
b H_{i}(\ol{b}) H_{j}(\ol{b})
+
\int_{a}^{\ol{b}}
H_{i}(v) H_{j}(v) \, dv
\nonumber \\ 
& & 
- 
a^{2} H_{i}(a) H_{j}(a) 
-
a H_{i}(a) 
\int_{a}^{\ol{b}}
H_{j}(v) \, dv
-
a b H_{i}(a) H_{j}(\ol{b}) 
\nonumber \\
& & 
- 
a H_{j}(a) 
\int_{a}^{\ol{b}}
H_{j}(v) \, dv
- 
\lp 
\int_{a}^{\ol{b}}
H_{i}(w) \, dw 
\rp
\lp 
\int_{a}^{\ol{b}}
H_{j}(v) \, dv
\rp 
-
b H_{j}(\ol{b}) 
\int_{a}^{\ol{b}}
H_{i}(v) \, dv 
\nonumber \\ 
& & 
- 
a b  H_{i}(\ol{b}) H_{j}(a) 
-
b H_{i}(\ol{b}) 
\int_{a}^{\ol{b}}
H_{j}(v) \, dv
-
b^{2} H_{i}(\ol{b}) H_{j}(\ol{b}) 
\nonumber \\ 
\therefore \ 
V_{11} 
& = & 
a H_{i}(a) H_{j}(a) 
+
b H_{i}(\ol{b}) H_{j}(\ol{b})
+
\int_{a}^{\ol{b}}
H_{i}(v) H_{j}(v) \, dv
- 
\Delta_{i} \, \Delta_{j},
\label{eqn:cgj17}
\nonumber 
\end{eqnarray}
as desired.
\end{proof}
\end{cor}

\begin{note}
The calculations of $c_{k}^{*}$,
for $1 \leq k \leq 3$, 
in Appendix A of \cite{MR2497558} 
are simply the special cases of 
Corollary \ref{cor:MTM1}.
\qed 
\end{note}

\section{Method of Winsorized Moments}
\label{sec:MWM} 

With the motivation of assigning special 
weights for a finite number of sample 
percentiles, 
the $L$-statistics given by
Eq. \eqref{eqn:S1} can be modified as 
\begin{align}
\label{eqn:S2}
\widehat{\mu}_{j} 
& :=
\dfrac{1}{n}
\sum_{i=1}^{n}
J_{j} \left( \dfrac{i}{n+1} \right) 
h_{j}(X_{i:n})
+ 
\sum_{i=1}^{m} 
c_{jn}^{(i)} 
h_{j} 
\left( 
X_{\left\lfloor n p_{j}^{(i)}\right\rfloor:n}
\right); 
\quad 
1 \le j \le k,
\end{align}
where 
$0 \le p_{j}^{(1)} < \cdots < p_{j}^{(m)} \le 1$
are the pre-specified percentiles and 
$c_{jn}^{(1)}, \ldots, c_{jn}^{(m)}$
are non-negative numbers such that
$
\displaystyle 
\lim_{n \to \infty} 
c_{jn}^{(i)} 
= 
c_{j}^{(i)}.
$
For example, 
we can select
\begin{align*}
c_{jn}^{(i)}
& \in
\left\{
\dfrac{\left\lfloor n \, p_{j}^{(i)}\right\rfloor}{n},
\dfrac{\left\lfloor n \, \overline{p}_{j}^{(i)}\right\rfloor}{n}
\right\},
\ \text{yielding} \ 
\displaystyle 
\lim_{n \to \infty} 
c_{jn}^{(i)} 
= 
c_{j}^{(i)}
\in 
\left\{
p_{j}^{(i)}, 
\overline{p}_{j}^{(i)}
\right\}, 
\
1 \le j \le k, \ 
1 \le i \le m.
\end{align*}
Define 
\begin{align}
\label{eqn:P3}
\mu_{j}
& = 
\int_{0}^{1}
{J_{j}(u)H_{j}(u) \, du
+ 
\sum_{i=1}^{m} c_{j}^{(i)} 
h_{j}
\left( 
F^{-1}
\left( 
p_{j}^{(i)}
\right)
\right); 
\quad 
1\leq{j}\leq{k}}.
\end{align}

\citet[][{\sc{Corollary}} 3]{MR0203874}
proved that the $k$-variate vector
$\sqrt{n}(\widehat{\bm\mu}-\bm{\mu})$,
converges in distribution to the 
$k$-variate normal random vector
with mean $\mathbf{0}$ and 
the variance-covariance matrix $\mathbf{\Sigma}:= \left[\sigma_{ij}^{2}\right]_{i,j=1}^{k}$ 
with the entries 
\begin{align}
\sigma_{ij}^{2} 
& = 
\int_{0}^{1} 
{\alpha_{i}(u)\alpha_{j}(u) \, du}  
\label{eqn:mwm_var_cov1}, 
\end{align}
where the function $\alpha$, 
like in Eq. \eqref{eqn:AlphaFun1},
is defined as
\begin{align}
\label{eqn:AlphaFun2}
\alpha_{j}(u)
& = 
\frac{1}{\ol{u}}
\left[ 
\int_{u}^{1}
\ol{v} \ g_{j}(v) \, dv
+ 
\sum_{i=1}^{m}
\ID \left\{{p}_{j}^{(i)} \ge u \right\}
c_{j}^{(i)} \ 
\overline{p}_{j}^{(i)} \ 
H' \left({p}_{j}^{(i)}\right)
\right];
\quad 
1\leq{j}\leq{k}.
\end{align}

Now, let $m = 2$, 
$p_{j}^{(1)} = a_{j}$, 
$p_{j}^{(2)} = \ol{b}_{j}$
for some 
$0 \le a_{j}, b_{j} \le 1$
with $ a_{j} < \ol{b}_{j}$.
Also, consider 
$
c_{jn}^{(1)}
=
\left\lfloor n a_{j} \right\rfloor/n
$
and 
$
c_{jn}^{(2)}
=
\left\lfloor n b_{j} \right\rfloor/n
$
yielding  
$
\displaystyle 
\lim_{n \to \infty} 
c_{jn}^{(1)} 
= 
c_{j}^{(1)}
=
a_{j}
$
and 
$
\displaystyle 
\lim_{n \to \infty} 
c_{jn}^{(2)} 
= 
c_{j}^{(2)}
=
b_{j}.
$
Define 
\begin{align}
\label{eqn:J_Fun2}
J_{j}(s) 
& =
\left\{ 
\begin{array}{ll}
1; & a_{j} < s < \ol{b}_{j}, \\
0; & \text{otherwise}. \\
\end{array}
\right. 
\quad 
1 \le j \le k. 
\end{align}

Eq. \eqref{eqn:S2}
is asymptotically equivalent to
the following equation,
see e.g., \citet[][p. 264]{MR595165}
\begin{align}
\widehat{\mu}_{j} 
& :=
\dfrac{1}{n}
\left[ 
\lf n a_{j} \rf 
h_{j}
\left( 
X_{\lf n a_{j} \rf +1:n}
\right) 
+
\sum_{i=\lf n a_{j} \rf +1}^{n-\lf n b_{j} \rf } 
h_{j}(X_{i:n}) 
+
\lf n b_{j} \rf  
h_{j}
\left(
X_{n - \lf n b_{j} \rf :n}) 
\right) 
\right].
\label{eqn:S3}
\end{align}

Similarly,
with the weights-generating function
given by Eq. \eqref{eqn:J_Fun2},
Eq. \eqref{eqn:P3},
takes the form 
\begin{align} 
\mu_{j}
& =
a_{j} 
h_{j}
\left(
F^{-1}(a_{j}) 
\right) 
+
\int_{a_{j}}^{\overline{b}_{j}} 
h_{j} \left( F^{-1}(v) \right)  \, dv 
+
b_{j} 
h_{j}
\left(
F^{-1}(\overline{b}_{j} )
\right).
\label{eqn:P4}
\end{align}

The expressions given by 
Eq. \eqref{eqn:S3} and 
Eq. \eqref{eqn:P4} are, 
respectively,
called sample- and population-
{\em winsorized moments}, 
see e.g., \citet[][p. 264, Example B (ii)]{MR595165}. 
Now, we match the sample and population
winsorized moments
from \eqref{eqn:S3} and \eqref{eqn:P4} 
to get the following system of equations for
$\theta_{1},\theta_{2},...,\theta_{k}$:
\begin{equation} 
\label{eq:match_mwm} 
\left\{
\begin{array}{lcl}
\mu_1 (\theta_1, \ldots, \theta_k) 
& = & 
\widehat{\mu}_1 \\
& \vdots & \\
\mu_k (\theta_1, \ldots, \theta_k) 
& = & 
\widehat{\mu}_k  \\
\end{array} \right.
\end{equation}

A solution, say 
$\widehat{\bm{\theta}}
=
\lp
\widehat{\theta}_{1},\widehat{\theta}_{2},
...,
\widehat{\theta}_{k}
\rp$, 
if it exists, 
to the system of equations (\ref{eq:match_mwm}) 
is called the {\textit{method of winsorized moments (MWM)}} 
estimator of $\bm{\theta}$. 
Thus, $\widehat{\theta}_{j}
=:
\tau_{j}
\left(
\widehat{\mu}_{1},\widehat{\mu}_{2},..., \widehat{\mu}_{k}
\right)$, $1\leq{j}\leq{k}$, 
are the MWM estimators of 
$\theta_{1},\theta_{2},...,\theta_{k}$.

With $m = 2$,
function $\alpha$, 
given by Eq. \eqref{eqn:AlphaFun2}
now takes the form
\begin{align}
\alpha_{j}(u)
& = 
\frac{1}{\ol{u}}
\left[ 
\int_{u}^{1}
\ol{v} \ g_{j}(v) \, dv
+ 
\sum_{i=1}^{2}
\ID \left\{{p}_{j}^{(i)} \ge u \right\}
c_{j}^{(i)} \ 
\overline{p}_{j}^{(i)} \ 
H_{j}' \left({p}_{j}^{(i)}\right)
\right] 
\nonumber \\
& = 
\frac{1}{\ol{u}}
\left[ 
\underbrace{
\int_{u}^{1}
\ol{v} \ 
g_{j}(v) \, dv}_{=: \psi_{1}(j)}
+ 
\underbrace{
\ID \left\{a_{j} \ge u \right\}
a_{j} \ 
\overline{a}_{j} \ 
H_{j}' \left(a_{j}\right)}_{=: \psi_{2}(j)}
+
\underbrace{
\ID \left\{\overline{b}_{j} \ge u \right\}
b_{j}^{2} \ 
H_{j}' \left(\overline{b}_{j}\right)}_{=: \psi_{3}(j)}
\right]
\nonumber \\[5pt]
& = 
\frac{
\psi_{1}(j)
+ 
\psi_{2}(j)
+
\psi_{3}(j)
}{\ol{u}}.
\label{eqn:AlphaFun5}
\end{align}

Then, 
from Eqs. \eqref{eqn:mwm_var_cov1} and \eqref{eqn:AlphaFun5},
it follows that
\begin{eqnarray}
\sigma_{ij}^{2} 
& = &
\int_{0}^{1} 
{\alpha_{i}(u)\alpha_{j}(u) \, du} 
\nonumber \\ 
& = &
\int_{0}^{1} 
\left( 
\dfrac{1}{\ol{u}}
\right)^{2} 
\left( 
\psi_{1}(i)
+ 
\psi_{2}(i)
+
\psi_{3}(i)
\right)
\left( 
\psi_{1}(j)
+ 
\psi_{2}(j)
+
\psi_{3}(j)
\right) \, du 
\nonumber \\
& =: & 
V_{11} + V_{12} + V_{13} 
+
V_{21} + V_{22} + V_{23}
+ 
V_{31} + V_{32} + V_{33},
\label{eqn:mwm_var_cov2}
\end{eqnarray}
where $V_{st}$, with $1 \leq s, t \leq 3$, 
are defined and evaluated in more detail,
which makes them more applicable and useful 
for different underlying distributions 
having at most a countable number 
of points of discontinuity.
\begin{enumerate}
\item 
Evaluating $V_{11}$:
\begin{eqnarray*}
V_{11}
& : = & 
\int_{0}^{1} 
\left( 
\dfrac{1}{\ol{u}}
\right)^{2} 
\psi_{1}(i) \psi_{1}(j) \, du 
\nonumber \\ 
& = & 
\int_{0}^{1} 
\left( 
\dfrac{1}{\ol{u}}
\right)^{2} 
\left[ 
\int_{u}^{1}
\ol{v} \ g_{i}(v) \, dv
\times 
\int_{u}^{1}
\ol{w} \ g_{j}(w) \, dw
\right] \, du 
\nonumber \\[5pt]
& = & 
\int_{a_{i}}^{\ol{b}_{i}} 
\int_{a_{j}}^{\ol{b}_{j}} 
\left[ 
v \wedge w - vw
\right]
H_{j}^{'}(v)
H_{i}^{'}(w) \, dv \, dw,
\end{eqnarray*}
which is already evaluated in 
Theorem \ref{thm:MTM_Var1}.

\item 
Now, we evaluate $V_{12}$.
Define 
$v \vee w := \max\{v,w\}$,
then we have
\begin{eqnarray}
V_{12} 
& : = & 
\int_{0}^{1} 
\left( 
\dfrac{1}{\ol{u}}
\right)^{2} 
\psi_{1}(i) \psi_{2}(j) \, du 
\nonumber \\
& = & 
\int_{0}^{1} 
\left( 
\dfrac{1}{\ol{u}}
\right)^{2}  
\ID \left\{a_{j} \ge u \right\}
a_{j} \ 
\overline{a}_{j} \ 
H_{j}' \left(a_{j}\right)
\left( 
\int_{u}^{1}
\ol{v} g_{i}(v) \, dv
\right) \, du 
\nonumber \\ 
& = & 
a_{j} \ 
\overline{a}_{j} \
H_{j}' \left(a_{j}\right)
\int_{0}^{1} 
\left( 
\dfrac{1}{\ol{u}}
\right)^{2}  
\ID \left\{a_{j} \ge u \right\}
\left( 
\int_{u}^{1}
\ol{v} \ g_{i}(v) \, dv
\right) \, du 
\nonumber \\
& = &
a_{j} \ 
\overline{a}_{j} \
H_{j}' \left(a_{j}\right)
\int_{0}^{1} 
\left( 
\dfrac{1}{\ol{u}}
\right)^{2}  
\ID \left\{a_{j} \ge u \right\}
\left( 
\int_{u}^{1}
J_{i}(v)H_{i}'(v) \, \ol{v} \, dv
\right) \, du 
\nonumber \\
& = &
a_{j} \ 
\overline{a}_{j} \
H_{j}' \left(a_{j}\right)
\int_{0}^{1} 
\left( 
\dfrac{1}{\ol{u}}
\right)^{2}  
\ID \left\{a_{j} \ge u \right\}
\left( 
\int_{a_{i} \vee u}^{\ol{b}_{i}}
H_{i}'(v) \, \ol{v} \, dv
\right) \, du
\nonumber \\
& = &
a_{j} \ 
\overline{a}_{j} \
H_{j}' \left(a_{j}\right)
\int_{0}^{1} 
\int_{a_{i} \vee u}^{\ol{b}_{i}}
\ID \left\{a_{j} \ge u \right\}
\left( 
\dfrac{1}{\ol{u}}
\right)^{2} 
H_{i}'(v) \, \ol{v} \, dv \, du
\nonumber \\
& = &
a_{j} \ 
\overline{a}_{j} \
H_{j}' \left(a_{j}\right)
\underbrace{
\int_{0}^{a_{j}} 
\int_{a_{i} \vee u}^{\ol{b}_{i}}
\left( 
\dfrac{1}{\ol{u}}
\right)^{2} 
H_{i}'(v) \, \ol{v} \, dv \, du}_{=: I_{12}}  
\nonumber \\
& = &
a_{j} \ 
\overline{a}_{j} \
H_{j}' \left(a_{j}\right)
I_{12}.
\label{eqn:V12_3}
\end{eqnarray} 

We evaluate \(I_{12}\) 
considering the inequality \eqref{eqn:abCondition1}.
\begin{eqnarray}
I_{12} 
& = &
\int_{0}^{a_{j}} 
\int_{a_{i} \vee u}^{\ol{b}_{i}}
\left( 
\dfrac{1}{\ol{u}}
\right)^{2}
H_{i}'(v) \, \ol{v} \, dv \, du  
\nonumber \\
& = & 
\int_{a_{j}}^{\ol{b}_{i}}
\int_{0}^{a_{j}}
\left( 
\dfrac{1}{\ol{u}}
\right)^{2}
H_{i}'(v) \, \ol{v} \, du \, dv 
\nonumber \\
& & 
+ 
\int_{a_{i}}^{a_{j}}
\int_{0}^{v}
\left( 
\dfrac{1}{\ol{u}}
\right)^{2}
H_{i}'(v) \, \ol{v} \, du \, dv 
\nonumber \\
& = &
\int_{a_{j}}^{\ol{b}_{i}}
H_{i}'(v) \, \ol{v}
\left(
\int_{0}^{a_{j}}
\left(
\dfrac{1}{\ol{u}}
\right)^{2} \, du 
\right) \, dv 
\nonumber \\
& & 
+ 
\int_{a_{i}}^{a_{j}}
H_{i}'(v) \, \ol{v}
\left(
\int_{0}^{v}
\left(
\dfrac{1}{\ol{u}}
\right)^{2} \, du 
\right) \, dv 
\nonumber \\
& = & 
\dfrac{a_{j}}{\ol{a}_{j}}
\int_{a_{j}}^{\ol{b}_{i}}
H_{i}'(v) \, \ol{v} \, dv 
+ 
\int_{a_{i}}^{a_{j}}
H_{i}'(v) \, \ol{v}
\left(
\dfrac{v}{\ol{v}} 
\right) \, dv 
\nonumber \\
& = & 
\dfrac{a_{j}}{\ol{a}_{j}}
\int_{a_{j}}^{\ol{b}_{i}}
H_{i}'(v) \, \ol{v} \, dv 
+ 
\int_{a_{i}}^{a_{j}}
v H_{i}'(v) \, dv 
\nonumber \\
& = & 
\dfrac{a_{j}}{\ol{a}_{j}}
\left[ 
b_{i} H_{i}(\ol{b}_{i}) 
-
\ol{a}_{j} H_{i}(a_{j})
+ 
\int_{a_{j}}^{\ol{b}_{i}} H_{i}(v) \, dv
\right] 
\nonumber \\
& & 
+ 
a_{j} H_{i}(a_{j}) 
-
a_{i} H_{i}(a_{i})
- 
\int_{a_{i}}^{a_{j}} H_{i}(v) \, dv.
\label{eqn:I12_2}
\end{eqnarray}

Thus, 
from Eq. \eqref{eqn:V12_3}, 
we have
\begin{eqnarray}
V_{12}
& = &
a_{j} \ 
\overline{a}_{j} \
H_{j}' \left(a_{j}\right)
I_{12} 
\nonumber \\ 
& = &
a_{j} \ 
\overline{a}_{j} \
H_{j}' \left(a_{j}\right)
\left\{
\dfrac{a_{j}}{\ol{a}_{j}}
\left[ 
b_{i} H_{i}(\ol{b}_{i}) 
-
\ol{a}_{j} H_{i}(a_{j})
+ 
\int_{a_{j}}^{\ol{b}_{i}} H_{i}(v) \, dv
\right]
\right. 
\nonumber \\
&  &
\left. 
+ 
a_{j} H_{i}(a_{j}) 
-
a_{i} H_{i}(a_{i})
- 
\int_{a_{i}}^{a_{j}} H_{i}(v) \, dv
\right\} 
\nonumber \\ 
& = & 
a_{j}^{2} \
H_{j}' \left(a_{j}\right)
\left[ 
b_{i} H_{i}(\ol{b}_{i}) 
-
\ol{a}_{j} H_{i}(a_{j})
+ 
\int_{a_{j}}^{\ol{b}_{i}} H_{i}(v) \, dv
\right] 
\nonumber \\
&  &
+
a_{j} \ 
\overline{a}_{j} \
H_{j}' \left(a_{j}\right)
\left\{ 
a_{j} H_{i}(a_{j}) 
-
a_{i} H_{i}(a_{i})
- 
\int_{a_{i}}^{a_{j}} H_{i}(v) \, dv
\right\}.
\nonumber 
\end{eqnarray}

\item 
Using a similar approach to that used for \( V_{12}\), 
we now calculate \( V_{13} \).
\begin{eqnarray}
V_{13} 
& : = & 
\int_{0}^{1} 
\left( 
\dfrac{1}{\ol{u}}
\right)^{2} 
\psi_{1}(i) \psi_{3}(j) \, du 
\nonumber \\
& = & 
\int_{0}^{1} 
\left( 
\dfrac{1}{\ol{u}}
\right)^{2}  
\ID \left\{\overline{b}_{j} \ge u \right\}
b_{j}^{2} \  
H_{j}' \left(\overline{b}_{j}\right)
\left( 
\int_{u}^{1}
\ol{v} \ g_{i}(v) \, dv
\right) \, du  
\nonumber \\
& = & 
b_{j}^{2} \ 
H_{j}' \left(\overline{b}_{j}\right)
\int_{0}^{1} 
\left( 
\dfrac{1}{\ol{u}}
\right)^{2}  
\ID \left\{\overline{b}_{j} \ge u \right\}
\left( 
\int_{u}^{1}
\ol{v} \ g_{i}(v) \, dv
\right) \, du  
\nonumber \\
& = & 
b_{j}^{2} \  
H_{j}' \left(\overline{b}_{j}\right)
\int_{0}^{1} 
\int_{a_{i} \vee u}^{\ol{b}_{i}}
\ID \left\{\overline{b}_{j} \ge u \right\}
\left( 
\dfrac{1}{\ol{u}}
\right)^{2} 
H_{i}'(v) \, \ol{v} \, dv \, du
\nonumber \\
& = &
b_{j}^{2} \ 
H_{j}' \left(\overline{b}_{j}\right)
\underbrace{
\int_{0}^{\overline{b}_{j}} 
\int_{a_{i} \vee u}^{\ol{b}_{i}}
\left( 
\dfrac{1}{\ol{u}}
\right)^{2} 
H_{i}'(v) \, \ol{v} \, dv \, du}_{=: I_{13}}  
\nonumber \\
& = & 
b_{j}^{2} \ 
H_{j}' \left(\overline{b}_{j}\right) 
I_{13}.
\label{eqn:V13_2} 
\end{eqnarray}

Again, considering inequality \eqref{eqn:abCondition1},
then similarly as in Eq. \eqref{eqn:I12_2}, 
it follows that 
\begin{eqnarray}
I_{13} 
& = &
\int_{0}^{\overline{b}_{j}} 
\int_{a_{i} \vee u}^{\overline{b}_{i}}
\left( 
\dfrac{1}{\ol{u}}
\right)^{2} 
H_{i}'(v) \, \ol{v} \, dv \, du 
\nonumber \\
& = & 
\dfrac{\overline{b}_{j}}{b_{j}}
\left[ 
b_{i} H_{i}(\ol{b}_{i}) 
-
b_{j} H_{i}(\overline{b}_{j})
+ 
\int_{\overline{b}_{j}}^{\ol{b}_{i}} H_{i}(v) \, dv
\right] 
\nonumber \\
& & 
+ 
\overline{b}_{j} H_{i}(\overline{b}_{j}) 
-
a_{i} H_{i}(a_{i})
- 
\int_{a_{i}}^{\overline{b}_{j}} H_{i}(v) \, dv.
\nonumber 
\end{eqnarray}

Thus, 
from Eq. \eqref{eqn:V13_2}, 
we have
\begin{eqnarray}
V_{13}
& = &
b_{j}^{2} \ 
H_{j}' \left(\overline{b}_{j}\right) 
I_{13} 
\nonumber \\ 
& = &
b_{j}^{2} \ 
H_{j}' \left(\overline{b}_{j}\right)
\left\{ 
\dfrac{\overline{b}_{j}}{b_{j}}
\left[ 
b_{i} H_{i}(\ol{b}_{i}) 
-
b_{j} H_{i}(\overline{b}_{j})
+ 
\int_{\overline{b}_{j}}^{\ol{b}_{i}} H_{i}(v) \, dv
\right] 
\right. 
\nonumber \\
& & 
\left. 
+ 
\overline{b}_{j} H_{i}(\overline{b}_{j}) 
-
a_{i} H_{i}(a_{i})
- 
\int_{a_{i}}^{\overline{b}_{j}} H_{i}(v) \, dv
\right\} 
\nonumber \\ 
& = &  
b_{j} \ \overline{b}_{j} \ 
H_{j}' \left(\overline{b}_{j}\right)
\left[ 
b_{i} H_{i}(\ol{b}_{i}) 
-
b_{j} H_{i}(\overline{b}_{j})
+ 
\int_{\overline{b}_{j}}^{\ol{b}_{i}} H_{i}(v) \, dv
\right] 
\nonumber \\
& & 
+
b_{j}^{2} \
H_{j}' \left(\overline{b}_{j}\right)
\left\{
\overline{b}_{j} H_{i}(\overline{b}_{j}) 
-
a_{i} H_{i}(a_{i})
- 
\int_{a_{i}}^{\overline{b}_{j}} H_{i}(v) \, dv
\right\}.
\label{eqn:V13_4}
\end{eqnarray}

\item 
Now, we evaluate $V_{21}$.
\begin{eqnarray}
V_{21} 
& : = & 
\int_{0}^{1} 
\left( 
\dfrac{1}{\ol{u}}
\right)^{2} 
\psi_{2}(i) \psi_{1}(j) \, du 
\nonumber \\
& = & 
\int_{0}^{1} 
\left( 
\dfrac{1}{\ol{u}}
\right)^{2}  
\ID \left\{a_{i} \ge u \right\}
a_{i} \ 
\overline{a}_{i} \ 
H_{i}' \left(a_{i}\right)
\left( 
\int_{u}^{1}
\ol{v} \ g_{j}(v) \, dv
\right) \, du 
\nonumber \\ 
& = & 
a_{i} \ 
\overline{a}_{i}
H_{i}' \left(a_{i}\right)
\int_{0}^{1} 
\left( 
\dfrac{1}{\ol{u}}
\right)^{2}  
\ID \left\{a_{i} \ge u \right\}
\left( 
\int_{u}^{1}
\ol{v} \ g_{j}(v) \, dv
\right) \, du
\nonumber \\
& = & 
a_{i} \ 
\overline{a}_{i} \
H_{i}' \left(a_{i}\right)
\int_{0}^{1} 
\left( 
\dfrac{1}{\ol{u}}
\right)^{2}  
\ID \left\{a_{i} \ge u \right\}
\left( 
\int_{u}^{1}
\ol{v} \, 
J_{j}(v)H_{j}'(v) \, dv
\right) \, du 
\nonumber \\
& = &
a_{i} \ 
\overline{a}_{i} \
H_{i}' \left(a_{i}\right)
\int_{0}^{1} 
\left( 
\dfrac{1}{\ol{u}}
\right)^{2}  
\ID \left\{a_{i} \ge u \right\}
\left( 
\int_{a_{j} \vee u}^{\ol{b}_{j}}
\ol{v} \, 
H_{j}'(v) \, dv
\right) \, du
\nonumber \\
& = &
a_{i} \ 
\overline{a}_{i} \
H_{i}' \left(a_{i}\right)
\int_{0}^{1} 
\int_{a_{j} \vee u}^{\ol{b}_{j}}
\ID \left\{a_{i} \ge u \right\}
\left( 
\dfrac{1}{\ol{u}}
\right)^{2} 
\ol{v} \, 
H_{j}'(v) \, dv \, du
\nonumber \\
& = &
a_{i} \ 
\overline{a}_{i} \
H_{i}' \left(a_{i}\right)
\underbrace{
\int_{0}^{a_{i}} 
\int_{a_{j} \vee u}^{\ol{b}_{j}}
\left( 
\dfrac{1}{\ol{u}}
\right)^{2} 
\ol{v} \, 
H_{j}'(v) \, dv \, du}_{=: I_{21}}  
\nonumber \\
& = &
a_{i} \ 
\overline{a}_{i} \
H_{i}' \left(a_{i}\right)
I_{21}.
\label{eqn:V21_2}
\end{eqnarray}

Noting the inequality \eqref{eqn:abCondition1},
i.e., $a_{i} \le a_{j}$,
we have 
\begin{align}
I_{21} 
& = 
\int_{0}^{a_{i}} 
\int_{a_{j} \vee u}^{\ol{b}_{j}}
\left( 
\dfrac{1}{\ol{u}}
\right)^{2} 
\ol{v} \, 
H_{j}'(v) \, dv \, du
\nonumber \\ 
& =
\int_{a_{j}}^{\ol{b}_{j}}
\ol{v} \, 
H_{j}'(v)
\left(
\int_{0}^{a_{i}}
\left( 
\dfrac{1}{\ol{u}}
\right)^{2} \, du
\right) \, dv 
\nonumber \\
& = 
\dfrac{a_{i}}{\ol{a}_{i}}
\int_{a_{j}}^{\ol{b}_{j}}
\ol{v} \, 
H_{j}'(v) \, dv 
\nonumber \\
& = 
\dfrac{a_{i}}{\ol{a}_{i}}
\left[ 
b_{j} H_{j}(\ol{b}_{j}) 
-
\ol{a}_{j} H_{j}(a_{j})
+ 
\int_{a_{j}}^{\ol{b}_{j}} H_{j}(v) \, dv
\right].
\label{eqn:I21_1}
\end{align}

Therefore, 
from Eqs. \eqref{eqn:V21_2} and \eqref{eqn:I21_1}, 
it follows that 
\begin{eqnarray*}
V_{21} 
& = & 
a_{i}^{2} \ 
H_{i}' \left(a_{i}\right)
\left[ 
b_{j} H_{j}(\ol{b}_{j}) 
-
\ol{a}_{j} H_{j}(a_{j})
+ 
\int_{a_{j}}^{\ol{b}_{j}} H_{j}(v) \, dv
\right].
\end{eqnarray*}

\item 
Next, we evaluate \(V_{22}\).
\begin{eqnarray}
V_{22} 
& : = & 
\int_{0}^{1} 
\left( 
\dfrac{1}{\ol{u}}
\right)^{2} 
\psi_{2}(i) \psi_{2}(j) \, du 
\nonumber \\
& = & 
\int_{0}^{1} 
\left( 
\dfrac{1}{\ol{u}}
\right)^{2}  
\ID \left\{a_{i} \ge u \right\}
a_{i} \ 
\overline{a}_{i} \ 
H_{i}' \left(a_{i}\right)
\ID \left\{a_{j} \ge u \right\}
a_{j} \ 
\overline{a}_{j} \ 
H_{j}' \left(a_{j}\right)
\nonumber \\ 
& = & 
a_{i} \ a_{j} \ 
\overline{a}_{i} \ \overline{a}_{j} \
H_{i}' \left(a_{i}\right) \  H_{j}' \left(a_{j}\right) 
\int_{0}^{1} 
\left( 
\dfrac{1}{\ol{u}}
\right)^{2}  
\ID \left\{a_{i} \ge u \right\}
\ID \left\{a_{j} \ge u \right\} \, du 
\nonumber \\ 
& = & 
a_{i} \ a_{j} \ 
\overline{a}_{i} \ \overline{a}_{j} \
H_{i}' \left(a_{i}\right) \  H_{j}' \left(a_{j}\right)
\left( 
\dfrac{1}{\ol{a}_{i} \wedge a_{j}} - 1
\right) 
\nonumber \\ 
& = & 
a_{i} \ a_{j} \ 
\overline{a}_{i} \ \overline{a}_{j} \
H_{i}' \left(a_{i}\right) \  H_{j}' \left(a_{j}\right) 
\begin{cases}
a_{i} \left(\overline{a}_{i}\right)^{-1}; & 
\mbox{if } a_{i} \le a_{j}, \\ 
a_{j} \left(\overline{a}_{j}\right)^{-1}; & 
\mbox{if } a_{i} > a_{j}
\end{cases} 
\nonumber \\ 
& = & 
\begin{cases}
a_{i}^{2} \ a_{j} \ \overline{a}_{j} \ 
H_{i}' \left(a_{i}\right) 
H_{j}' \left(a_{j}\right); & 
\mbox{if } a_{i} \le a_{j}, \\ 
a_{i} \ a_{j}^{2} \ \overline{a}_{i} \ 
H_{i}' \left(a_{i}\right) 
H_{j}' \left(a_{j}\right); & 
\mbox{if } a_{i} > a_{j}.
\end{cases} 
\nonumber
\end{eqnarray}

\item 
Next, we evaluate \(V_{23}\).
\begin{eqnarray}
V_{23} 
& : = & 
\int_{0}^{1} 
\left( 
\dfrac{1}{\ol{u}}
\right)^{2} 
\psi_{2}(i) \psi_{3}(j) \, du 
\nonumber \\
& = & 
\int_{0}^{1} 
\left( 
\dfrac{1}{\ol{u}}
\right)^{2}  
\ID \left\{a_{i} \ge u \right\}
a_{i} \ 
\overline{a}_{i} \ 
H_{i}' \left(a_{i}\right)
\ID \left\{a_{j} \ge u \right\}
\ID \left\{\overline{b}_{j} \ge u \right\}
b_{j}^{2} \ 
H_{j}' \left(\overline{b}_{j}\right) \, du
\nonumber \\
& = & 
a_{i} \ \overline{a}_{i} \ 
b_{j}^{2} \ 
H_{i}' \left(a_{i}\right) 
H_{j}' \left(\overline{b}_{j}\right) 
\int_{0}^{1} 
\left( 
\dfrac{1}{\ol{u}}
\right)^{2}  
\ID \left\{a_{i} \ge u \right\}
\ID \left\{\overline{b}_{j} \ge u \right\} \, du 
\nonumber \\
& = & 
a_{i} \ \overline{a}_{i} \ 
b_{j}^{2} \ 
H_{i}' \left(a_{i}\right) 
H_{j}' \left(\overline{b}_{j}\right) 
\int_{0}^{a_{i} \wedge \overline{b}_{j}} 
\left( 
\dfrac{1}{\ol{u}}
\right)^{2} \, du 
\nonumber \\
& = & 
a_{i} \ \overline{a}_{i} \ 
b_{j}^{2} \ 
H_{i}' \left(a_{i}\right) 
H_{j}' \left(\overline{b}_{j}\right)  
\left.
\left(
\dfrac{1}{\ol{u}}
\right)
\right|_{0}^{a_{i} \wedge \overline{b}_{j}}  
\nonumber \\
& = & 
a_{i} \ \overline{a}_{i} \ 
b_{j}^{2} \ 
H_{i}' \left(a_{i}\right)
H_{j}' \left(\overline{b}_{j}\right)  
\left(
\dfrac{a_{i} \wedge \overline{b}_{j}}
{1 - a_{i} \wedge \overline{b}_{j}}
\right).
\label{eqn:V23_1}
\end{eqnarray}

\item 
The evaluation of \(V_{31}\) is similar to that of \(V_{13}\). 
\begin{eqnarray}
V_{31} 
& : = & 
\int_{0}^{1} 
\left( 
\dfrac{1}{\ol{u}}
\right)^{2} 
\psi_{3}(i) \psi_{1}(j) \, du 
\nonumber \\
& = & 
\int_{0}^{1} 
\left( 
\dfrac{1}{\ol{u}}
\right)^{2} 
\ID \left\{\overline{b}_{i} \ge u \right\}
b_{i}^{2} \ 
H_{i}' \left(\overline{b}_{i}\right)
\left( 
\int_{u}^{1}
\ol{v} \ g_{j}(v) \, dv
\right) \, du
\nonumber \\
& = & 
b_{i}^{2} \ 
H_{i}' \left(\overline{b}_{i}\right)
\int_{0}^{1} 
\left( 
\dfrac{1}{\ol{u}}
\right)^{2} 
\ID \left\{\overline{b}_{i} \ge u \right\}
\left( 
\int_{u}^{1}
\ol{v} \ g_{j}(v) \, dv
\right) \, du.
\nonumber 
\end{eqnarray}

Noting $a_{j} \le \overline{b}_{i}$
from inequality \eqref{eqn:abCondition1}, 
then following Eq. \eqref{eqn:V13_4},
we have
\begin{eqnarray*}
V_{31}
& = &  
b_{i} \ \overline{b}_{i} \ 
H_{i}' \left(\overline{b}_{i}\right)
\left[ 
b_{j} H_{j}(\ol{b}_{j}) 
-
b_{i} H_{j}(\overline{b}_{i})
+ 
\int_{\overline{b}_{i}}^{\ol{b}_{j}} H_{j}(v) \, dv
\right] 
\nonumber \\
& & 
+
b_{i}^{2} \
H_{i}' \left(\overline{b}_{i}\right)
\left\{
\overline{b}_{i} H_{j}(\overline{b}_{i}) 
-
a_{j} H_{j}(a_{j})
- 
\int_{a_{j}}^{\overline{b}_{i}} H_{j}(v) \, dv
\right\}.
\end{eqnarray*}

\item 
Next,
we evaluate \(V_{32}\), 
which closely resembles the calculation of 
\(V_{23}\) as outlined in Eq. \eqref{eqn:V23_1}.
Thus, we have
\begin{eqnarray}
V_{32} 
& := &
\int_{0}^{1} 
\left( 
\dfrac{1}{\ol{u}}
\right)^{2} 
\psi_{3}(i) \psi_{2}(j) \, du 
\nonumber \\
& = & 
\int_{0}^{1} 
\left( 
\dfrac{1}{\ol{u}}
\right)^{2}  
\ID \left\{\overline{b}_{i} \ge u \right\}
b_{i}^{2} \ 
H_{i}' \left(\overline{b}_{i}\right)
\ID \left\{a_{j} \ge u \right\}
\ID \left\{a_{j} \ge u \right\}
a_{j} \ 
\overline{a}_{j} \ 
H_{j}' \left(a_{j}\right) 
\nonumber \\
& = & 
a_{j} \ \overline{a}_{j} \ 
b_{i}^{2} \ 
H_{j}' \left(a_{j}\right) \ H_{i}' \left(\overline{b}_{i}\right) 
\int_{0}^{1} 
\left( 
\dfrac{1}{\ol{u}}
\right)^{2}  
\ID \left\{a_{j} \ge u \right\}
\ID \left\{\overline{b}_{i} \ge u \right\} \, du 
\nonumber \\
& = & 
a_{j} \ \overline{a}_{j} \ 
b_{i}^{2} \ 
H_{j}' \left(a_{j}\right) \ H_{i}' \left(\overline{b}_{i}\right)  
\left(
\dfrac{a_{j} \wedge \overline{b}_{i}}
{1 - a_{j} \wedge \overline{b}_{i}}
\right).
\nonumber 
\end{eqnarray}

\item 
Finally, we evaluate \(V_{33}\),
which differs from the other quantities 
and is presented as follows.
\begin{eqnarray}
V_{33} 
& : = & 
\int_{0}^{1} 
\left( 
\dfrac{1}{\ol{u}}
\right)^{2} 
\psi_{3}(i) \psi_{3}(j) \, du 
\nonumber \\ 
& = & 
\int_{0}^{1} 
\left( 
\dfrac{1}{\ol{u}}
\right)^{2} 
\ID \left\{\overline{b}_{i} \ge u \right\}
b_{i}^{2} \ 
H_{i}' \left(\overline{b}_{i}\right)
\ID \left\{\overline{b}_{j} \ge u \right\}
b_{j}^{2} \ 
H_{j}' \left(\overline{b}_{j}\right) \, du  
\nonumber \\ 
& = & 
b_{i}^{2} \ 
b_{j}^{2} \
H_{i}' \left(\overline{b}_{i}\right) \ 
H_{j}' \left(\overline{b}_{j}\right)
\int_{0}^{1} 
\left( 
\dfrac{1}{\ol{u}}
\right)^{2} 
\ID \left\{\overline{b}_{i} \ge u \right\}
\ID \left\{\overline{b}_{j} \ge u \right\} \, du 
\nonumber \\
& = & 
b_{i}^{2} \ 
b_{j}^{2} \
H_{i}' \left(\overline{b}_{i}\right) \ 
H_{j}' \left(\overline{b}_{j}\right) 
\left( 
\dfrac{1}{1 - \overline{b}_{i} \wedge \overline{b}_{j}} - 1
\right) 
\nonumber \\ 
& = & 
b_{i}^{2} \ 
b_{j}^{2} \
H_{i}' \left(\overline{b}_{i}\right) \ 
H_{j}' \left(\overline{b}_{j}\right)
\begin{cases}
\overline{b}_{i} \ b_{i}^{-1}; & 
\mbox{if } b_{i} \ge b_{j}, \\ 
\overline{b}_{j} \ b_{j}^{-1}; & 
\mbox{if } b_{i} < b_{j}
\end{cases} 
\nonumber \\ 
& = & 
\begin{cases}
b_{i} \ \overline{b}_{i} \ b_{j}^{2} \
H_{i}' \left(\overline{b}_{i}\right)  
H_{j}' \left(\overline{b}_{j}\right); & 
\mbox{if } b_{i} \ge b_{j}, \\ 
b_{i}^{2} \ b_{j} \ \overline{b}_{j} \
H_{i}' \left(\overline{b}_{i}\right)  
H_{j}' \left(\overline{b}_{j}\right); & 
\mbox{if } b_{i} < b_{j}.
\end{cases} 
\nonumber 
\end{eqnarray}
\end{enumerate} 

\begin{cor}
\label{cor:EstablishA1}
Given the special condition specified by 
inequality \eqref{eqn:abSituation1}, 
prove that 
\begin{align*}
V_{11}
& = 
a H_{i}(a) H_{j}(a) 
+
b H_{i}(\ol{b}) H_{j}(\ol{b})
+
\int_{a}^{\ol{b}}
H_{i}(v) H_{j}(v) \, dv
- 
\Delta_{i} \, \Delta_{j},
\nonumber \\
V_{12}  
& = 
a^{2} \
H_{j}' \left(a\right)
\left[ 
b H_{i}(\ol{b}) 
-
\ol{a} H_{i}(a)
+ 
\int_{a}^{\ol{b}} H_{i}(v) \, dv
\right],
\nonumber \\
V_{13} 
& = 
b^{2} 
H_{j}' \left(\overline{b}\right)
\left[
\overline{b} H_{i}(\overline{b}) 
-
a H_{i}(a)
- 
\int_{a}^{\overline{b}} H_{i}(v) \, dv
\right],
\nonumber \\
V_{21} 
& = 
a^{2} 
H_{i}' \left(a\right)
\left[ 
b H_{j}(\ol{b}) 
-
\ol{a} H_{j}(a)
+ 
\int_{a}^{\ol{b}} H_{j}(v) \, dv
\right],
\nonumber \\
V_{22} 
& = 
a^{3} \ \overline{a} \ 
H_{i}' \left( a \right) 
H_{j}' \left( a \right),
\quad 
V_{23} 
=
a^{2}   
b^{2} 
H_{i}' \left(a\right)
H_{j}'\left( \overline{b} \right), 
\nonumber \\ 
V_{31} 
& = 
b^{2} 
H_{i}' \left(\overline{b}\right)
\left[
\overline{b} H_{j}(\overline{b}) 
-
a H_{j}(a)
- 
\int_{a}^{\overline{b}} H_{j}(v) \, dv
\right],
\nonumber \\
V_{32} 
& = 
a^{2} 
b^{2} 
H_{j}' \left(a\right) 
H_{i}' \left(\overline{b}\right), 
\quad 
V_{33} 
= 
b^{3} \ \overline{b}
H_{i}' \left(\overline{b}\right)  
H_{j}' \left(\overline{b}\right).
\end{align*}

\begin{proof}
A proof of $V_{11}$ is already provided 
in Corollary \ref{cor:MTM1}. 
The proofs of the remaining quantities follow 
directly from their corresponding general cases.
\end{proof}
\end{cor}

\begin{note}
The calculations of 
\(\widehat{A}_{i,j}^{(k)}\) for \(1 \leq k \leq 4\), 
as presented in Lemma A.1 from \cite{MR3758788}, 
can simply be reexpressed using the
results of Corollary \ref{cor:EstablishA1}.
\qed 
\end{note}

\section{Concluding Remarks}
\label{sec:Conclusion}

This paper establishes computational formulas 
for the asymptotic variance of two robust 
$L$-estimators: the method of trimmed moments (MTM) and 
the method of winsorized moments (MWM). 
We demonstrate that two asymptotic approaches for 
MTM are equivalent with a specific choice 
of the weights-generating function, 
Theorem \ref{thm:MTM_Theorem1}.
This enhances the usability of these estimators 
across various statistical settings,
particularly in complex actuarial scenarios
such as payment-per-payment and payment-per-loss. 
Regardless of the nature of the observed sample data,
such as truncation or censorship,
if the cumulative distribution function 
of the underlying distribution is definable, 
then for any trimming/winsorizing proportions,
the asymptotic computational formulas;
Theorem \ref{thm:MTM_Var1}, 
Corollaries \ref{cor:MTM1} and \ref{cor:EstablishA1}, 
Eqs. \eqref{eqn:cgj9} and \eqref{eqn:mwm_var_cov2},
can be applied.
Moreover, 
since distortion risk measure can be expressed
in the form of Eq. \eqref{eqn:P1}, 
as discussed in \cite{bpv23, MR3039556}, 
the computational formulas developed in this paper
are also applicable to establishing the 
asymptotic normality properties of distortion risk measure.


\newpage 

\begin{thebibliography}{2024}
\baselineskip 4.80mm
\setlength{\bibsep}{5.0pt plus 0.00ex}
\bibitem[\protect\astroncite{Bernard et~al.}{2023}]{bpv23}
Bernard, C., Pesenti, S.~M., and Vanduffel, S. (2023).
\newblock Robust distortion risk measures.
\newblock {\em Mathematical Finance}.

\bibitem[\protect\astroncite{Brazauskas et~al.}{2007}]{MR2416885}
Brazauskas, V., Jones, B.~L., and Zitikis, R. (2007).
\newblock Robustification and performance evaluation of empirical risk measures and other vector-valued estimators.
\newblock {\em METRON--International Journal of Statistics}, 65(2):175--199.

\bibitem[\protect\astroncite{Brazauskas et~al.}{2009}]{MR2497558}
Brazauskas, V., Jones, B.~L., and Zitikis, R. (2009).
\newblock Robust fitting of claim severity distributions and the method of trimmed moments.
\newblock {\em Journal of Statistical Planning and Inference}, 139(6):2028--2043.

\bibitem[\protect\astroncite{Brazauskas and Kleefeld}{2009}]{MR2591318}
Brazauskas, V. and Kleefeld, A. (2009).
\newblock Robust and efficient fitting of the generalized {P}areto distribution with actuarial applications in view.
\newblock {\em Insurance: Mathematics \& Economics}, 45(3):424--435.

\bibitem[\protect\astroncite{Chernoff et~al.}{1967}]{MR0203874}
Chernoff, H., Gastwirth, J.~L., and Johns, Jr., M.~V. (1967).
\newblock Asymptotic distribution of linear combinations of functions of order statistics with applications to estimation.
\newblock {\em Annals of Mathematical Statistics}, 38(1):52--72.

\bibitem[\protect\astroncite{Dhaene et~al.}{2012}]{MR3039556}
Dhaene, J., Kukush, A., Linders, D., and Tang, Q. (2012).
\newblock Remarks on quantiles and distortion risk measures.
\newblock {\em European Actuarial Journal}, 2(2):319--328.

\bibitem[\protect\astroncite{Folland}{1999}]{MR1681462}
Folland, G.~B. (1999).
\newblock {\em Real Analysis: Modern Techniques and Their Applications}.
\newblock John Wiley \& Sons, New York, second edition.

\bibitem[\protect\astroncite{Hosking}{1990}]{MR1049304}
Hosking, J. R.~M. (1990).
\newblock {$L$}-moments: analysis and estimation of distributions using linear combinations of order statistics.
\newblock {\em Journal of the Royal Statistical Society. Series B. Methodological}, 52(1):105--124.

\bibitem[\protect\astroncite{Klugman et~al.}{2019}]{MR3890025}
Klugman, S.~A., Panjer, H.~H., and Willmot, G.~E. (2019).
\newblock {\em Loss Models: From Data to Decisions}.
\newblock John Wiley \& Sons, Hoboken, NJ, fifth edition.

\bibitem[\protect\astroncite{Poudyal}{2021}]{MR4263275}
Poudyal, C. (2021).
\newblock Robust estimation of loss models for lognormal insurance payment severity data.
\newblock {\em ASTIN Bulletin -- The Journal of the International Actuarial Association}, 51(2):475--507.

\bibitem[\protect\astroncite{Poudyal and Brazauskas}{2022}]{pb22}
Poudyal, C. and Brazauskas, V. (2022).
\newblock Robust estimation of loss models for truncated and censored severity data.
\newblock {\em Variance}, 15(2):1--20.

\bibitem[\protect\astroncite{Serfling}{1980}]{MR595165}
Serfling, R.~J. (1980).
\newblock {\em Approximation Theorems of Mathematical Statistics}.
\newblock John Wiley \& Sons, New York.

\bibitem[\protect\astroncite{Zhao et~al.}{2018}]{MR3758788}
Zhao, Q., Brazauskas, V., and Ghorai, J. (2018).
\newblock Robust and efficient fitting of severity models and the method of {W}insorized moments.
\newblock {\em ASTIN Bulletin -- The Journal of the International Actuarial Association}, 48(1):275--309.
\end{thebibliography}
\end{document}